\documentclass{amsart}
\usepackage{amsmath, amssymb}
\theoremstyle{plain}
 \newtheorem{thm}{Theorem}[section]
 \newtheorem{lem}[thm]{Lemma}
 \newtheorem{cor}[thm]{Corollary}
 \newtheorem{prop}[thm]{Proposition}

\theoremstyle{definition}
  \newtheorem{defn}{Definition}[section]
  
  \newtheorem{assumption}{Assumption}[section]

\theoremstyle{remark}
  \newtheorem{rem}{Remark}[section]

\newcommand{\ci}[2]{\cite[#1]{#2}}

\begin{document}
\title[Measure equivalence rigidity of the mapping class group]{Measure equivalence rigidity \\ of the mapping class group}
\author{Yoshikata Kida}
\address{Max-Planck-Institut f\"ur Mathematik, Vivatsgasse 7, Bonn 53111, Germany}
\email{kida@mpim-bonn.mpg.de}

\begin{abstract}
We show that the mapping class group of a compact orientable surface with higher complexity has the following extreme rigidity in the sense of measure equivalence: if the mapping class group is measure equivalent to a discrete group, then they are commensurable up to finite kernel. Moreover, we describe all lattice embeddings of the mapping class group into a locally compact second countable group. We also obtain similar results for finite direct products of mapping class groups.
\end{abstract}

\maketitle

\footnote[0]{{\it Date}: July 19, 2006.}

\section{Introduction}\label{sec-int}
The purpose of this paper is to establish a new rigidity theorem for the mapping class group in terms of measure equivalence among discrete groups. In this paper, by a discrete group we mean a discrete and countable one. Measure equivalence was introduced by Gromov \cite{gromov2} as follows:

\begin{defn}
We say that two discrete groups $\Gamma$ and $\Lambda$ are measure equivalent if there exists a measure-preserving action of $\Gamma \times \Lambda$ on a standard Borel space $(\Sigma, m)$ with a $\sigma$-finite positive measure such that both restricted actions to $\Gamma$ and $\Lambda$ are essentially free and have a fundamental domain of finite measure. The space $(\Sigma, m)$ (equipped with the $\Gamma \times \Lambda$-action) is then called a ME coupling of $\Gamma$ and $\Lambda$. 
\end{defn}

It is easy to see that measure equivalence defines an equivalence relation among discrete groups. One typical example of two measure equivalent groups is given by any two lattices in the same locally compact second countable (lcsc) group, which is the main geometric motivation for the introduction of this notion. Commensurability up to finite kernel is the equivalence relation for discrete groups defined by declaring two groups in an exact sequence $1\rightarrow A\rightarrow B\rightarrow C\rightarrow 1$ of discrete groups to be equivalent if the third group is finite. Any two discrete groups which are commensurable up to finite kernel are measure equivalent.

Measure equivalence between two groups has another equivalent formulation in terms of orbit equivalence (see Subsection \ref{subsec-ME oe}), which has been studied for a long time and is closely related to ergodic theory and the theory of von Neumann algebras. The first magnificent result about orbit equivalence is due to Ornstein-Weiss \cite{ow} following Dye \cite{dye}, \cite{dye2} and it can be stated in terms of measure equivalence as follows: a discrete group is measure equivalent to $\mathbb{Z}$ if and only if it is an infinite amenable group. (This result was generalized for amenable discrete measured equivalence relations by Connes-Feldman-Weiss \cite{cfw}.) Note that Zimmer \cite{zim2} extended the superrigidity theorem for semisimple Lie groups of non-compact type due to Mostow-Margulis to one in the context of orbit equivalence, which is called his cocycle superrigidity theorem, and he classified lattices in simple Lie groups of real rank at least $2$ up to measure equivalence. In addition to this classification, his cocycle superrigidity theorem has many applications to various rigidity phenomena of higher rank lattices.

Recently, the study of measure equivalence and orbit equivalence is very rapidly developing. Furman's beautiful rigidity result in \cite{furman1} completely determines the class of discrete groups measure equivalent to some higher rank lattices. Namely, if a discrete group $\Lambda$ is measure equivalent to a lattice in a connected simple Lie group $G$ with finite center and real rank at least $2$, then there exists a homomorphism from $\Lambda$ onto a lattice in ${\rm Ad}G$ whose kernel is finite. Gaboriau's discovery in \cite{gab-l2} that $\ell^{2}$-Betti numbers for discrete groups are invariant under measure equivalence in a certain sense leads to surprising progress in the classification problem of measure equivalence because these numerical invariants are defined for all discrete groups and are computable for various discrete groups arising geometrically. Popa showed in \cite{popa-mal2}, \cite{popa} that the Bernoulli shifts of Kazhdan groups have various strong rigidity properties in terms of orbit equivalence. It is remarkable that he treated all Kazhdan groups, which is a very large class of groups. The reader should be referred to \cite{gab-survey}, \cite{shalom-survey}, \cite{vaes} for more details of recent development of measure equivalence and orbit equivalence.

Let $M=M_{g, p}$ be a connected compact orientable surface of type $(g, p)$, that is, of genus $g$ and with $p$ boundary components. Throughout the paper, a surface is assumed to be connected, compact and orientable unless otherwise stated. The mapping class group $\Gamma(M)$ of $M$ is defined as the group of isotopy classes of all orientation-preserving diffeomorphisms of $M$. The extended mapping class group $\Gamma(M)^{\diamond}$ of $M$ is the group of isotopy classes of all diffeomorphisms of $M$, which contains $\Gamma(M)$ as a subgroup of index $2$. Let $\kappa(M)=3g+p-4$ be the complexity of $M$. If $\kappa(M)>0$, we say that $M$ has higher complexity. Let $C=C(M)$ be the curve complex for a surface $M$. In \cite{kida}, we obtain some classification result of $\Gamma(M)$ in terms of measure equivalence and give various examples of discrete groups not measure equivalent to $\Gamma(M)$. In the proof, we established fundamental methods to study subrelations in a discrete measured equivalence relation arising from a standard action of $\Gamma(M)$, where a standard action means an essentially free, measure-preserving action on a standard Borel space with a finite positive measure. The curve complex played one of the most important roles in the study of them. Using these methods, we show the following rigidity theorem for $\Gamma(M)$, which can be viewed as an analogue of Furman's rigidity theorem and is the main result of this paper. Let ${\rm Aut}(C)$ denote the automorphism group of the simplicial complex $C$. Note that ${\rm Aut}(C)$ and $\Gamma(M)$ are commensurable up to finite kernel (see Theorem \ref{thm-cc-auto}).

\begin{thm}\label{thm-main}
If a discrete group $\Lambda$ is measure equivalent to the mapping class group $\Gamma(M)$ with $\kappa(M)>0$, then there exists a homomorphism $\rho \colon \Lambda \rightarrow {\rm Aut}(C)$ whose kernel and cokernel are both finite. 
\end{thm}

This theorem completely determines the class of discrete groups measure equivalent to $\Gamma(M)$ and provides the first example of infinite discrete groups with such an extreme rigidity in the theory of measure equivalence. Remark that uniform and non-uniform lattices in Lie groups treated in Furman's rigidity result are not commensurable up to finite kernel because they are not quasi-isometric \cite{farb}. Hence, if $\Gamma$ is a lattice as in Furman's rigidity theorem, then there exists a discrete group which is measure equivalent to $\Gamma$ and is not commensurable up to finite kernel with $\Gamma$. 

\begin{rem}
Although both Furman's and Popa's rigidity theorems are concerned with discrete groups satisfying (or related to) Kazhdan's property, the mapping class group for a surface of genus at most $2$ does not have Kazhdan's property or more strongly, it contains a subgroup of finite index which admits a quotient isomorphic to a non-abelian free group of finite rank (see \cite{kork}). It is unknown whether the mapping class groups for other surfaces have Kazhdan's property or not. 
\end{rem}

Theorem \ref{thm-main} completes the classification of mapping class groups up to measure equivalence.

\begin{thm}\label{cor-classification}
Suppose that $M^{1}$, $M^{2}$ are distinct surfaces of type $(g_{1}, p_{1})$, $(g_{2}, p_{2})$, respectively with $\kappa(M^{1}), \kappa(M^{2})>0$ and $g_{1}\leq g_{2}$. Moreover, assume that $\Gamma(M^{1})$ and $\Gamma(M^{2})$ are measure equivalent. Then we only have the following two possibilities: $((g_{1}, p_{1}), (g_{2}, p_{2}))=((0, 5), (1, 2)), ((0, 6), (2, 0))$. 
\end{thm}

\begin{rem}
Note that if $\kappa(M)<0$ and $M\neq M_{1, 0}$, then $\Gamma(M)$ is finite. Both $\Gamma(M_{1, 0})$ and $\Gamma(M_{1, 1})$ are isomorphic to $SL(2, \mathbb{Z})$ and $\Gamma(M_{0, 4})$ is commensurable up to finite kernel with $SL(2, \mathbb{Z})$. It is known that $\Gamma(M_{0, 5})$ and $\Gamma(M_{1, 2})$ (resp. $\Gamma(M_{0, 6})$ and $\Gamma(M_{2, 0}$)) are commensurable up to finite kernel (see Theorem \ref{thm-cc-auto}).  
\end{rem}

Let $\Gamma$ be a lattice in a connected simple Lie group $G$ with finite center and real rank at least $2$. In the proof of Furman's rigidity theorem in \cite{furman1}, the main ingredient is to prove the following (see \cite[Theorem 4.1]{furman1}): let $(\Omega, \omega)$ be a self ME coupling of $\Gamma$ (i.e., a ME coupling of $\Gamma$ and $\Gamma$). Then there exists an essentially unique almost $\Gamma \times \Gamma$-equivariant Borel map $\Psi \colon \Omega \rightarrow {\rm Aut}({\rm Ad}G)$, which means that
\[\Psi((\gamma, \gamma')x)={\rm Ad}(\gamma)\Psi(x){\rm Ad}(\gamma')^{-1}\]
for any $\gamma, \gamma' \in \Gamma$ and a.e.\ $x\in \Omega$. Furman substantially used Zimmer's cocycle superrigidity theorem for the construction of $\Psi$. On the other hand, we can show that for a self ME coupling $(\Sigma, m)$ of $\Gamma(M)$, there exists an essentially unique almost $\Gamma(M) \times \Gamma(M)$-equivariant Borel map $\Phi \colon \Sigma \rightarrow {\rm Aut}(C)$, which means that
\[\Phi((\gamma, \gamma')x)=\pi(\gamma)\Phi(x)\pi(\gamma')^{-1}\]
for any $\gamma, \gamma'\in \Gamma(M)$ and a.e.\ $x\in \Sigma$. Here, $\pi \colon \Gamma(M)^{\diamond}\rightarrow {\rm Aut}(C)$ is the natural homomorphism. This construction of $\Phi$ is the heart of the proof of Theorem \ref{thm-main}. In Section \ref{sec-gr-gen}, we give an outline of the construction of $\Phi$.

After the construction of the map $\Phi$, we apply Furman's and Monod-Shalom's techniques in \cite{furman1}, \cite{ms} for higher rank lattices and non-trivial direct products of finitely many groups with special properties (e.g., word-hyperbolic ones), which are applicable to a more general situation. More precisely, given a ME coupling $(\Sigma', m')$ of $\Gamma(M)$ and a discrete group $\Lambda$, we construct the self ME coupling $\Sigma'\times_{\Lambda}\Lambda \times_{\Lambda}\check{\Sigma}'$. Using their techniques for the equivariant Borel map from this self ME coupling into ${\rm Aut}(C)$, one can find a homomorphism $\rho$ as in Theorem \ref{thm-main}.

Moreover, we consider the same problem as above for a finite direct product $\Gamma(M_{1})\times \cdots \times \Gamma(M_{n})$ of mapping class groups $\Gamma(M_{i})$ with $\kappa(M_{i})>0$ for all $i$. Monod-Shalom introduced in \cite{ms} the class $\mathcal{C}$ consisting of discrete groups $\Gamma$ which admit a mixing unitary representation $\pi$ on a Hilbert space such that the second bounded cohomology $H_{b}^{2}(\Gamma, \pi)$ of $\Gamma$ with coefficient $\pi$ does not vanish. They showed in that paper that a non-trivial finite direct product of discrete groups in $\mathcal{C}$ satisfies various measurable rigidity properties. The class $\mathcal{C}$ contains a large number of discrete groups arising geometrically (e.g., word-hyperbolic ones) and whether a discrete group is in $\mathcal{C}$ or not is invariant under measure equivalence. Hamenst\"adt proved in \cite{ham2} that the mapping class group of a surface with higher complexity is contained in $\mathcal{C}$, so we can obtain various measurable rigidity theorems as in \cite{ms} for direct products of mapping class groups. 

Following Monod-Shalom's ingenious technique treating fundamental domains for actions on ME couplings, we can find an essentially unique almost equivariant Borel map from a self ME coupling of a direct product of $\Gamma(M_{i})$ into the  direct product of ${\rm Aut}(C(M_{i}))$. In the same way as above, we can show the following theorem similar to Theorem \ref{thm-main} for a direct product of $\Gamma(M_{i})$: 

\begin{thm}\label{thm-main-prod}
Let $n$ be a positive integer and let $M_{i}$ be a surface with $\kappa(M_{i})>0$ for $i\in \{ 1, \ldots, n\}$. If a discrete group $\Lambda$ is measure equivalent to the direct product $\Gamma(M_{1})\times \cdots \times \Gamma(M_{n})$, then there exists a homomorphism $\Lambda \rightarrow {\rm Aut}(C(M_{1}))\times \cdots \times {\rm Aut}(C(M_{n}))$ whose kernel and cokernel are both finite.
\end{thm}

In \cite{furman3}, Furman gave another application of the map $\Psi$ mentioned above. He gave an explicit description of a lcsc group containing a lattice isomorphic to a lattice in a simple Lie group of higher rank. Roughly speaking, he showed that such a lcsc group can be built only from the ambient Lie group or from the lattice itself and their actions on a compact group. Following Furman, we describe a lcsc group containing a lattice isomorphic to the mapping class group as follows. We fix notations as follows: let $n\in \mathbb{N}$ and let $M_{i}$ be a surface with  $\kappa(M_{i})>0$ for $i\in \{ 1, \ldots, n\}$. Put 
\[G={\rm Aut}(C(M_{1}))\times \cdots \times {\rm Aut}(C(M_{n}))\]
and let $G_{i}$ be equal to $\Gamma(M_{i})^{\diamond}$ or ${\rm Aut}(C(M_{i}))$. Put $G_{0}=G_{1}\times \cdots \times G_{n}$ and let $\pi \colon G_{0}\rightarrow G$ be the natural homomorphism.

\begin{thm}\label{thm-int-lattice-emb}
Suppose that $\Gamma$ is a subgroup of finite index in $G_{0}$ and there is a lattice embedding $\sigma \colon \Gamma \rightarrow H$ into a lcsc group $H$, that is, $\sigma$ is an injective homomorphism such that $\sigma(\Gamma)$ is a lattice in $H$. Then we have a subgroup $H_{0}$ of finite index in $H$ containing $\sigma(\Gamma)$ and a compact normal subgroup $K$ of $H_{0}$ satisfying the following:
\begin{enumerate}
\item[(i)] $[H: H_{0}]\leq [G: \pi(\Gamma)]$;
\item[(ii)] The action of $H$ on $K$ by conjugation induces via $\sigma$ an action of $\Gamma$ on $K$ and hence a semi-direct product $\Gamma \ltimes K$. Let $p\colon \Gamma \ltimes K\rightarrow H_{0}$ be the homomorphism defined by $\Gamma \ni \gamma \mapsto \sigma(\gamma)$ and $K\ni k\mapsto k$. Then $p$ is surjective;
\item[(iii)] For $\gamma \in \Gamma$ and $k\in K$, we have $(\gamma, k)\in \ker(p)$ if and only if $\pi(\gamma)=e$ and $k=\sigma(\gamma)^{-1}$. In particular, if the kernel of the restriction of $\pi$ to $\Gamma$ is trivial, then $p$ is an isomorphism.
\end{enumerate}
\end{thm} 

This theorem says that there exist no interesting lattice embeddings of the mapping class group. The following corollary can easily be shown:

\begin{cor}\label{cor-int-lattice-emb}
Let $\Gamma$ be a subgroup of finite index in $G_{0}$ and suppose that we have a lattice embedding of $\Gamma$ into a lcsc group $H$. Then the image of $\Gamma$ is cocompact in $H$ and $H$ has infinitely many connected components.
\end{cor}

It follows from this corollary that any subgroup of finite index of the mapping class group for a surface with higher complexity can not be isomorphic to a lattice in a semisimple Lie group, which was proved first by Kaimanovich-Masur \cite{kai-mas}. They showed more generally that any sufficiently large subgroup of the mapping class group can not be isomorphic to a lattice in a semisimple Lie group. In this direction, Farb-Masur \cite{farb-mas}, Bestvina-Fujiwara \cite{bestvina-fujiwara} and Yeung \cite{yeung} studied homomorphisms from a lattice in some semisimple Lie group into the mapping class group and concluded that their images are finite.

In a subsequent paper \cite{kida-OE}, we give an application of the existence of an equivariant Borel map from a self ME coupling of the mapping class group, following Furman's idea in \cite{furman2}. In \cite{kida-OE}, we establish orbit equivalence rigidity of ergodic standard actions of the mapping class group and give a new example of a discrete measured equivalence relation which can not arise from any standard action of a discrete group. Moreover, we give uncountably many explicit examples of ergodic standard actions of the mapping class group which are mutually non-orbit equivalent, using certain generalized Bernoulli shifts of the mapping class group.  

\vspace{1em}

\noindent {\bf Acknowledgement.} The author is grateful to Professor Ursula Hamenst\"adt for reading the first draft of this paper very carefully and giving many valuable suggestions. The author also thanks the Max Planck Institute for Mathematics at Bonn for its warm hospitality.


\section{Preliminaries}

\subsection{The mapping class group}\label{subsec-mcg}

In this subsection, we recall fundamental facts about the mapping class group and several geometric objects related to it. We refer the reader to \cite{FLP}, \cite{ivanov1}, \cite{ivanov2} or Sections 3.1, 3.2, 4.3 and 4.5 in \cite{kida} and the references therein for the material of this subsection.

Let $M=M_{g, p}$ be a surface of genus $g$ and with $p$ boundary components. Let $\Gamma(M)$, $\Gamma(M)^{\diamond}$ be the mapping class group and the extended one of $M$, respectively, defined as in Section \ref{sec-int}. Let $\kappa(M)=3g+p-4$ be the complexity of $M$. When $\kappa(M)>0$, we say that $M$ has higher complexity.

For a surface $M$, let $V(C)=V(C(M))$ be the set of all non-trivial isotopy classes of non-peripheral simple closed curves on $M$. Let $S(M)$ denote the set of all non-empty finite subsets of $V(C)$ which can be realized disjointly on $M$ at the same time. If $\kappa(M)>0$, we can define the curve complex $C=C(M)$ as a simplicial complex whose vertex set is $V(C)$ and simplex set is $S(M)$. Remark that when $\kappa(M)=0$, we can define the curve complex of $M$ in a slightly different way so that its vertex set $V(C)$ is given in the same way as above. If $\kappa(M)\geq 0$, then $\Gamma(M)^{\diamond}$ has the natural and simplicial action on $C$ and $C$ is connected and has infinite diameter. Moreover, when $C$ is equipped with a natural combinatorial metric, it is hyperbolic in the sense of Gromov (see \cite{masur-minsky1}).

Let $M$ be a surface with $\kappa(M)\geq 0$ and let $i\colon V(C)\times V(C)\rightarrow \mathbb{N}$ be the geometric intersection number. Let $\mathcal{MF}=\mathcal{MF}(M)$ be the space of measured foliations on $M$. Let $\mathcal{PMF}=\mathcal{PMF}(M)$ be the space of projective measured foliations on $M$, which is also called the Thurston boundary and is homeomorphic to the sphere of dimension $6g-7+2p$. Note that $S(M)$ can naturally be embedded into $\mathcal{PMF}$. The geometric intersection number $i$ can continuously and $\mathbb{R}_{>0}$-homogeneously be extended to a function $\mathcal{MF}\times \mathcal{MF}\rightarrow \mathbb{R}_{\geq 0}$ in the following sense:
\[i(r_{1}F_{1}, r_{2}F_{2})=r_{1}r_{2}i(F_{1}, F_{2})\]
for any $r_{1}, r_{2}\in \mathbb{R}_{>0}$ and $F_{1}, F_{2}\in \mathcal{MF}$. Hence, for two elements $F_{1}, F_{2}\in \mathcal{PMF}$, whether $i(F_{1}, F_{2})=0$ or $\neq 0$ makes sense. It is clear that $\Gamma(M)^{\diamond}$ acts continuously on both $\mathcal{MF}$ and $\mathcal{PMF}$ and 
\[i(gF_{1}, gF_{2})=i(F_{1}, F_{2})\]
for any $g\in \Gamma(M)^{\diamond}$ and $F_{1}, F_{2}\in \mathcal{MF}$ (or $\mathcal{PMF}$). Let 
\[\mathcal{MIN}=\{ F\in \mathcal{PMF}:i(F, \alpha)\neq 0\ {\rm for\ any\ }\alpha \in V(C)\}\]  
be the set of all minimal measured foliations on $M$, which is a $\Gamma(M)^{\diamond}$-invariant Borel subset of $\mathcal{PMF}$.

The Thurston boundary $\mathcal{PMF}$ is an ideal boundary of the Teichm\"uller space of $\mathcal{T}=\mathcal{T}(M)$ for $M$. The union $\overline{\mathcal{T}}=\mathcal{T}\cup \mathcal{PMF}$ is called the Thurston compactification of the Teichm\"uller space, which is homeomorphic to the closed Euclidean ball of dimension $6g-6+2p$ whose boundary corresponds to $\mathcal{PMF}$.

For $g\in \Gamma(M)$, let us denote by
\[{\rm Fix}(g)=\{ x\in \overline{\mathcal{T}}: gx=x\}\]
the fixed point set of $g$. Each element $g\in \Gamma(M)$ is classified as follows in terms of its fixed points on $\overline{\mathcal{T}}$ (see Exp\'ose 9, \S V, Th\'eor\`eme and Exp\'ose 11, \S 4, Th\'eor\`eme in \cite{FLP}):
\begin{enumerate}
\item[(i)] $g$ has finite order and has a fixed point on $\mathcal{T}$;
\item[(ii)] $g$ is pseudo-Anosov, which means that ${\rm Fix}(g)$ consists of exactly two points in $\mathcal{MIN}$;
\item[(iii)] $g$ has infinite order and is reducible, which means that there exists $\sigma \in S(M)$ such that $g\sigma =\sigma$.
\end{enumerate}
These three types are mutually exclusive. We say that $F\in \mathcal{PMF}$ is a pseudo-Anosov foliation if $F$ is a fixed point for some pseudo-Anosov element. The set of all pseudo-Anosov elements is dense in $\mathcal{PMF}$. 

Since the curve complex $C$ is hyperbolic, we can consider its boundary $\partial C$ at infinity, which is not compact. There exists a natural $\Gamma(M)$-equivariant continuous map $\mathcal{MIN}\rightarrow \partial C$, which is injective on the set of all uniquely ergodic measured foliations. This set contains all pseudo-Anosov foliations (see \cite{ham}, \cite{kla} or \cite[Section 3.2]{kida}). 

A pseudo-Anosov element $g\in \Gamma(M)$ has the following remarkable dynamics on $\overline{\mathcal{T}}$ (see \cite[Theorem 7.3.A]{ivanov2}): the two fixed points $F_{\pm}(g)\in \mathcal{MIN}$ of $g$ satisfy that if $U$ is any neighborhood of $F_{+}(g)$ in $\overline{\mathcal{T}}$ and $K$ is any compact set in $\overline{\mathcal{T}}\setminus \{ F_{-}(g)\}$, then $g^{n}(K)\subset U$ for all sufficiently large $n\in \mathbb{N}$. 

Using the above classification of elements of $\Gamma(M)$, McCarthy-Papadopoulos \cite{mc-pa} classified subgroups $\Gamma$ of $\Gamma(M)$ as follows:
\begin{enumerate}
\item[(i)] $\Gamma$ is finite;
\item[(ii)] there exists a pseudo-Anosov element $g\in \Gamma$ such that $h\{ F_{\pm}(g)\} =\{ F_{\pm}(g)\}$ for any $h\in \Gamma$. In this case, $\Gamma$ is virtually cyclic;
\item[(iii)] there exists $\sigma \in S(M)$ such that $g\sigma =\sigma$ for any $g\in \Gamma$. In this case, $\Gamma$ is said to be reducible;
\item[(iv)] $\Gamma$ contains an independent pair $\{ g_{1}, g_{2}\}$ of pseudo-Anosov elements, which means $\{ F_{\pm}(g_{1})\} \cap \{ F_{\pm}(g_{2})\} =\emptyset$. In this case, $\Gamma$ contains a non-abelian free subgroup and is said to be sufficiently large.
\end{enumerate}

We recall the canonical reduction system (CRS) for a subgroup of $\Gamma(M)$, which plays an important role in the study of reducible subgroups. We refer the reader to \cite[Chapter 7]{ivanov1} for more details about CRS's. For $\sigma \in S(M)$, we denote by $M_{\sigma}$ for simplicity the surface obtained by cutting along a realization of curves in $\sigma$. For an integer $m\geq 3$, let $\Gamma(M; m)$ be the subgroup of $\Gamma(M)$ consisting of all elements which act on the homology group $H_{1}(M; \mathbb{Z}/m\mathbb{Z})$ trivially. This subgroup has the following notable properties (see Theorem 1.2 and Corollaries 1.5, 1.8, 3.6 in \cite{ivanov1} or \cite[Section 4.3]{kida}): 

\begin{thm}\label{thm-pure}
In the above notation, the following assertions hold:
\begin{enumerate}
\item[(i)] $\Gamma(M;m)$ is a torsion-free subgroup of finite index in $\Gamma(M)$.
\item[(ii)] If $g\in \Gamma(M;m)$ and $F\in \mathcal{PMF}$ satisfy $g^{n}F=F$ for some $n\in \mathbb{Z}\setminus \{ 0\}$, then $gF=F$.
\item[(iii)] If $g\in \Gamma(M;m)$ and $\sigma \in S(M)$ satisfy $g^{n}\sigma =\sigma$ for some $n\in \mathbb{Z}\setminus \{ 0\}$, then $g\alpha =\alpha$ for any $\alpha \in S(M)$ and $g$ preserves each component of $M_{\sigma}$ and of the boundary of $M$.
\end{enumerate}
\end{thm}

When we consider the problem of measure equivalence in subsequent sections, we study actions of (a finite index subgroup of) $\Gamma(M; m)$ instead of $\Gamma(M)$.

Let $\Gamma$ be a subgroup of $\Gamma(M; m)$. A curve $\alpha \in V(C)$ is called an essential reduction class for $\Gamma$ if the following two conditions are satisfied:
\begin{enumerate}
\item[(i)] $g\alpha =\alpha$ for any $g\in \Gamma$;
\item[(ii)] if $\beta \in V(C)$ satisfies $i(\alpha, \beta)\neq 0$, then there exists $g\in \Gamma$ such that $g\beta \neq \beta$.
\end{enumerate}
The canonical reduction system (CRS) $\sigma(\Gamma)$ for $\Gamma$ is defined to be the set of all essential reduction classes for $\Gamma$, which is either an element in $S(M)$ or empty. We can define the CRS for a general subgroup $\Gamma$ of $\Gamma(M)$ as the CRS for $\Gamma \cap \Gamma(M;m)$, which is independent of $m$. It is known that an infinite subgroup $\Gamma$ of $\Gamma(M)$ is reducible if and only if $\sigma(\Gamma)$ is non-empty (see \cite[Corollary 7.17]{ivanov1}).

Given a subgroup $\Gamma$ of $\Gamma(M;m)$ and $\sigma \in S(M)$ with $g\sigma =\sigma$ for any $g\in \Gamma$, thanks to Theorem \ref{thm-pure} (iii), we can define a natural homomorphism
\[p_{\sigma}\colon \Gamma \rightarrow \prod_{Q}\Gamma(Q), \]
where $Q$ runs through all components of $M_{\sigma}$. 

\begin{lem}[\ci{Lemma 2.1 (1)}{blm}, \ci{Corollary 4.1.B}{ivanov2}]\label{lem-ker-dehn}
The kernel of $p_{\sigma}$ is contained in the subgroup of $\Gamma(M)$ generated by Dehn twists about all curves in $\sigma$. 
\end{lem}

For each component $Q$ of $M_{\sigma}$, let $p_{Q}\colon \Gamma \rightarrow \Gamma(Q)$ be the composition of $p_{\sigma}$ with the projection onto $\Gamma(Q)$. It is known that if a subgroup $\Gamma$ of $\Gamma(M;m)$ is reducible and $Q$ is a component of $M_{\sigma(\Gamma)}$, then the image $p_{Q}(\Gamma)$ either is trivial or contains a pseudo-Anosov element in $\Gamma(Q)$ (see \cite[Corollary 7.18]{ivanov1}). If $p_{Q}(\Gamma)$ is trivial, infinite amenable or non-amenable, then $Q$ is said to be T, IA or IN, respectively.

\begin{lem}\label{lem-red-crs-group}
Let $\Gamma$ be a finite index subgroup of $\Gamma(M)$ and define a subgroup
\[\Gamma_{\sigma}=\{ g\in \Gamma : g\sigma =\sigma \}\]
for $\sigma \in S(M)$. Then the CRS for $\Gamma_{\sigma}$ is equal to $\sigma$.
\end{lem}

This lemma follows easily from \cite[Theorem 7.16]{ivanov1} because any component of $M_{\sigma}$ which is not a pair of pants is IN for $\Gamma_{\sigma}$ if $\Gamma$ is a finite index subgroup of $\Gamma(M; m)$ with an integer $m\geq 3$.

\begin{lem}\label{lem-1-crs}
Let $\Gamma$ be an infinite subgroup of $\Gamma(M;m)$ and let $\alpha \in V(C(M))$. Assume that $g\alpha =\alpha$ for all $g\in \Gamma$. If for each component $Q$ of $M_{\alpha}$, we have $g\beta =\beta$ for any $\beta \in V(C(Q))$ and $g\in \Gamma$, then the CRS for $\Gamma$ is $\{ \alpha \}$.
\end{lem}

\begin{proof}
Since $\Gamma$ is infinite and reducible, the CRS $\sigma(\Gamma)$ for $\Gamma$ is non-empty. Let $\delta \in \sigma(\Gamma)$. We show that $\delta =\alpha$. Let $Q$ be a component of $M_{\alpha}$. If $\delta \in V(C(Q))$, then there exists $\beta \in V(C(Q))$ with $i(\beta, \delta)\neq 0$. By assumption, $\beta$ is invariant for $\Gamma$, which contradicts the assumption that $\delta$ is an essential reduction class for $\Gamma$. Thus, either $i(\delta, \alpha)\neq 0$ or $\delta =\alpha$. The former case can not happen because $\alpha$ is invariant for $\Gamma$ and $\delta \in \sigma(\Gamma)$.
\end{proof}


\subsection{The automorphism group of the curve complex}\label{subsec-aut-cc}

Let $M$ be a surface with $\kappa(M)>0$. Then we have the natural homomorphism $\pi \colon \Gamma(M)^{\diamond}\rightarrow {\rm Aut}(C)$. 

\begin{thm}[\cite{ivanov3}, \cite{luo}]\label{thm-cc-auto}
Let $M$ be a surface with $\kappa(M)>0$.
\begin{enumerate}
\item[(i)] If $M$ is neither $M_{1, 2}$ nor $M_{2, 0}$, then $\pi$ is an isomorphism.
\item[(ii)] If $M=M_{1, 2}$, then the image of $\pi$ is a subgroup of ${\rm Aut}(C)$ with index $5$ and $\ker(\pi)$ is the subgroup generated by a hyperelliptic involution, which is isomorphic to $\mathbb{Z}/2\mathbb{Z}$.  
\item[(iii)] If $M=M_{2, 0}$, then $\pi$ is surjective and $\ker(\pi)$ is the subgroup generated by a hyperelliptic involution, which is isomorphic to $\mathbb{Z}/2\mathbb{Z}$.
\item[(iv)] The two curve complexes $C(M_{0, 5})$, $C(M_{1, 2})$ (resp. $C(M_{0, 6})$, $C(M_{2, 0})$) are isomorphic as simplicial complexes.
\end{enumerate}
\end{thm}

\begin{thm}\label{thm-icc}
Let $\Gamma$ be a subgroup of finite index in ${\rm Aut}(C)$. For each $g_{0}\in {\rm Aut}(C)\setminus \{ e\}$, the set $\{ gg_{0}g^{-1}: g\in \Gamma \}$ consists of infinitely many elements.
\end{thm}

\begin{proof}
We may assume that $M\neq M_{1, 2}, M_{2, 0}$ by Theorem \ref{thm-cc-auto} (iv). By Theorem \ref{thm-cc-auto} (i), we identify ${\rm Aut}(C)$ with $\Gamma(M)^{\diamond}$. Let $g_{0}\in \Gamma(M)^{\diamond}$ and assume that the set $\{ gg_{0}g^{-1}: g\in \Gamma \}$ consists of only finitely many elements. Then note that for any infinite subset $\{ h_{n}\}_{n\in \mathbb{N}}$ of $\Gamma$, there exists an infinite subsequence $\{ n_{i}\}_{i\in \mathbb{N}}$ of $\mathbb{N}$ such that $h_{n_{i}}g_{0}h_{n_{i}}^{-1}=h_{n_{j}}g_{0}h_{n_{j}}^{-1}$ for each $i$, $j$. Put 
\[{\rm Fix}(g_{0})=\{ x\in \overline{\mathcal{T}}: g_{0}x=x\},\]
which is a non-empty closed subset of $\overline{\mathcal{T}}$. 

Assume ${\rm Fix}(g_{0})\not\supset \mathcal{PMF}$. If we deduce a contradiction, then the inclusion ${\rm Fix}(g_{0})\supset \mathcal{PMF}$ holds, which implies $g_{0}=e$ and completes the proof. Since the set of pseudo-Anosov foliations is dense in $\mathcal{PMF}$, we can find an independent pair $\{ g_{1}, g_{2}\}\subset \Gamma(M)$ of pseudo-Anosov elements such that $F_{+}(g_{1}), F_{+}(g_{2})\in \mathcal{PMF}\setminus {\rm Fix}(g_{0})$. Using the assumption that $\Gamma$ is a subgroup of finite index in $\Gamma(M)^{\diamond}$, we may assume that $g_{1}, g_{2}\in \Gamma$.

Let $s\in {\rm Fix}(g_{0})$. It follows from the above remark that there exists an infinite increasing subsequence $\{ n_{i}\}$ of $\mathbb{N}$ such that $g_{1}^{-n_{i}}g_{0}g_{1}^{n_{i}}=g_{0}$ for each $i$. Then we have $s=g_{0}s=g_{1}^{-n_{i}}g_{0}g_{1}^{n_{i}}s$, which implies $g_{1}^{n_{i}}s\in {\rm Fix}(g_{0})$ for each $i$. If $s\neq F_{-}(g_{1})$, then $g_{1}^{n_{i}}s\rightarrow F_{+}(g_{1})\in \mathcal{PMF}\setminus {\rm Fix}(g_{0})$ as $i\rightarrow \infty$, which is a contradiction. Thus, $s=F_{-}(g_{1})$. Similarly, we can show that if $s\in {\rm Fix}(g_{0})$, then $s=F_{-}(g_{2})$, which is a contradiction because the pair $\{ g_{1}, g_{2}\}$ is independent. 
\end{proof}


\subsection{Measure equivalence and orbit equivalence}\label{subsec-ME oe}

In this subsection, we recall the construction of weakly orbit equivalent actions from a ME coupling given in \cite[Section 3]{furman2}. We refer the reader to \cite{ana} for the terminology of an $r$-discrete measured groupoid and its amenability. We fix notations as follows: given an $r$-discrete measured groupoid $\mathcal{G}$ on a standard finite measure space $(X, \mu)$ (i.e., a standard Borel space with a finite positive measure) and a Borel subset $A\subset X$ with positive measure, we denote by
\[(\mathcal{G})_{A}=\{ \gamma \in \mathcal{G}: r(\gamma), s(\gamma)\in A\}\]
the restricted groupoid to $A$, where $r, s\colon \mathcal{G}\rightarrow X$ are the range, source maps, respectively. For $x, y\in X$, let 
\[\mathcal{G}^{x}_{y}=\{ \gamma \in \mathcal{G}: r(\gamma)=x, s(\gamma)=y\}\]
and let $e_{x}\in \mathcal{G}_{x}^{x}$ denote the unit. If $A$ is a Borel subset of $X$, then $\mathcal{G}A$ denotes the saturation defined by
\[\mathcal{G}A=\{ r(\gamma)\in X:  \gamma \in \mathcal{G}, s(\gamma)\in A\},\]
which is a Borel subset of $X$.

Let $(\Sigma, m)$ be a ME coupling of discrete groups $\Gamma$ and $\Lambda$, and choose fundamental domains $Y, X\subset \Sigma$ for the $\Gamma$-, $\Lambda$-actions. Remark that we have a natural $\Gamma$-action on $X$ equipped with the restricted measure of $m$ because $X$ can be identified with the quotient space $\Sigma /\Lambda$ as Borel spaces. Similarly, we have a natural $\Lambda$-action on $Y$. In order to distinguish from the original $\Gamma$-, $\Lambda$-actions on $\Sigma$, we denote the $\Gamma$-, $\Lambda$-actions on $X$, $Y$ by $\gamma \cdot x$, $\lambda \cdot y$, respectively, using a dot. Note that we can choose $X$ and $Y$ so that $A=X\cap Y$ satisfies $\Gamma \cdot A=X$ and $\Lambda \cdot A=Y$ up to null sets. In what follows, suppose that $X$ and $Y$ satisfy this condition.

Let $\mathcal{G}=\Gamma \ltimes X$ (resp. $\mathcal{H}=\Lambda \ltimes Y$) be the $r$-discrete measured groupoid on $(X, \mu)$ (resp. $(Y, \nu))$ constructed from the above action. We can define cocycles
\[\alpha \colon \Gamma \times X\rightarrow \Lambda,\ \ \beta \colon \Lambda \times Y\rightarrow \Gamma\]
so that $\gamma \cdot x=\alpha(\gamma, x)\gamma x\in X$ and $\lambda \cdot y=\beta(\lambda, y)\lambda y\in Y$ for any $\gamma \in \Gamma$, $\lambda \in \Lambda$ and a.e.\ $x\in X$, $y\in Y$. Let
\[p\colon X\rightarrow Y,\ \ q\colon Y\rightarrow X\]
be the Borel maps defined by
\[p(x)=\Gamma x\cap Y,\ \ q(y)=\Lambda y\cap X\]
for $x\in X$ and $y\in Y$. Note that both $p$ and $q$ are the identity on $A=X\cap Y$ and 
\[p(\gamma \cdot x)=\alpha(\gamma, x)\cdot p(x),\ \ q(\lambda \cdot y)=\beta(\lambda, y)\cdot q(y)\]
for any $\gamma \in \Gamma$, $\lambda \in \Lambda$ and a.e.\ $x\in X$, $y\in Y$. Define groupoid homomorphisms
\begin{align*}
f&\colon \mathcal{G}\ni (\gamma, x)\mapsto (\alpha(\gamma, x), p(x))\in \mathcal{H}\\ 
g&\colon \mathcal{H}\ni (\lambda, y)\mapsto (\beta(\lambda, y), q(y))\in \mathcal{G}.
\end{align*}
Note that $\beta(\alpha(\gamma, x), x)=\gamma$ for any $\gamma \in \Gamma$ and a.e.\  $x\in A$ with $\gamma \cdot x\in A$, and $\alpha(\beta(\lambda, y), y)=\lambda$ for any $\lambda \in \Lambda$ and a.e.\ $y\in A$ with $\lambda \cdot y\in A$. Therefore, we obtain the following:

\begin{prop}
The groupoid homomorphisms
\[f\colon (\mathcal{G})_{A}\rightarrow (\mathcal{H})_{A},\ \ g\colon (\mathcal{H})_{A}\rightarrow (\mathcal{G})_{A}\]
satisfy $g\circ f={\rm id}$ and $f\circ g={\rm id}$.
\end{prop}

This proposition implies that the two actions of $\Gamma$ on $X$ and $\Lambda$ on $Y$ are weakly orbit equivalent. 

Consider the $\Gamma \times \Lambda$-action on $X\times \Lambda$ defined by
\[\gamma \lambda(x, \lambda')=(\gamma \cdot x, \alpha(\gamma, x)\lambda'\lambda^{-1}).\]
It is easy to check the following lemma:

\begin{lem}\label{lem-ME iso}
The Borel map $\Sigma \rightarrow X\times \Lambda$ defined by $\lambda x\mapsto (x, \lambda^{-1})$ for $x\in X$ and $\lambda \in \Lambda$ is Borel isomorphic and $\Gamma \times \Lambda$-equivariant. 
\end{lem}

Conversely, we know the following theorem. For simplicity, a standard action of a discrete group means an essentially free, measure-preserving Borel action of it on a standard finite measure space.

\begin{thm}[\ci{Theorem 3.3}{furman2}]
If two discrete groups $\Gamma$ and $\Lambda$ have ergodic standard actions on $(X, \mu)$ and $(Y, \nu)$ which are weakly orbit equivalent, then we can construct a ME coupling $(\Sigma, m)$ of $\Gamma$ and $\Lambda$ such that the $\Gamma$-actions on $X$ and $\Lambda \backslash \Sigma$ (resp. the $\Lambda$-actions on $Y$ and $\Gamma \backslash \Sigma$) are isomorphic.
\end{thm}


\subsection{Normal subgroupoids} 

In this subsection, we introduce the notion of normal subgroupoids of an $r$-discrete measured groupoid, based on \cite{fsz}, \cite[Subsection 4.6.1]{kida}. This notion is a generalization of normal subrelations of a discrete measured equivalence relation and also a generalization of normal subgroups of a discrete group.

Let $\mathcal{G}$ be an $r$-discrete measured groupoid on a standard finite measure space $(X, \mu)$ and $r, s\colon \mathcal{G}\rightarrow X$ be the range, source maps, respectively. Let $\mathcal{S}$ be a subgroupoid of $\mathcal{G}$. In this paper, a subgroupoid of $\mathcal{G}$ means a Borel subgroupoid of $\mathcal{G}$ whose unit space is the same as the one for $\mathcal{G}$. Let us denote by ${\rm End}_{\mathcal{S}}(\mathcal{G})$ the set of all Borel maps $\varphi \colon {\rm dom}(\varphi)\rightarrow \mathcal{G}$ from a Borel subset of $X$ such that 
\begin{enumerate}
\item[(i)] $s(\varphi(x))=x$ for a.e.\ $x\in {\rm dom}(\varphi)$;
\item[(ii)] $\gamma \in \mathcal{S}$ if and only if $\varphi(r(\gamma))\gamma \varphi(s(\gamma))^{-1}\in \mathcal{S}$ for a.e.\ $\gamma \in (\mathcal{G})_{{\rm dom}(\varphi)}$.
\end{enumerate}
Let $[[\mathcal{G}]]_{\mathcal{S}}$ be the subset of ${\rm End}_{\mathcal{S}}(\mathcal{G})$ consisting of all $\varphi$ such that the map ${\rm dom}(\varphi)\ni x\mapsto r(\varphi(x))\in X$ is injective a.e.\ on the domain. When $\mathcal{S}$ is trivial, we write $[[\mathcal{G}]]$ instead of $[[\mathcal{G}]]_{\mathcal{S}}$. 

We define the composition $\psi \circ \varphi \colon {\rm dom}(\psi \circ \varphi)\rightarrow \mathcal{G}$ of two elements $\varphi, \psi \in {\rm End}_{\mathcal{S}}(\mathcal{G})$ by 
\[{\rm dom}(\psi \circ \varphi)=\{ x\in {\rm dom}(\varphi): r(\varphi(x))\in {\rm dom}(\psi)\},\] 
\[\psi \circ \varphi(x)=\psi(r(\varphi(x)))\varphi(x)\]
for $x\in {\rm dom}(\psi \circ \varphi)$. It is easy to check that $\psi \circ \varphi \in {\rm End}_{\mathcal{S}}(\mathcal{G})$.

\begin{defn}
A subgroupoid $\mathcal{S}$ of an $r$-discrete measured groupoid $\mathcal{G}$ on a standard finite measure space $(X, \mu)$ is said to be normal in $\mathcal{G}$ if the following condition is satisfied: there exists a countable family $\{ \phi_{n}\}$ of maps in ${\rm End}_{\mathcal{S}}(\mathcal{G})$ such that for a.e.\ $\gamma \in \mathcal{G}$, we can find $\phi_{n}$ in the family satisfying $r(\gamma)\in {\rm dom}(\phi_{n})$ and $\phi_{n}(r(\gamma))\gamma \in \mathcal{S}$. In this case, we write $\mathcal{S}\vartriangleleft \mathcal{G}$ and we call $\{ \phi_{n}\}$ a family of normal choice functions for the pair $(\mathcal{G}, \mathcal{S})$.
\end{defn}

The following two lemmas give natural examples of normal subgroupoids:

\begin{lem}
Suppose that a discrete group $G$ has a non-singular action on $(X, \mu)$ and $H$ is a normal subgroup of $G$. Let $\mathcal{G}_{G}$ and $\mathcal{G}_{H}$ be the groupoids generated by the actions of $G$ and $H$, respectively. Then the subgroupoid $\mathcal{G}_{H}$ is normal in $\mathcal{G}_{G}$.
\end{lem}

The following two lemmas can be proved by using \cite[Theorem 3.9]{kida}:

\begin{lem}\label{lem-iso-normal}
Let $\mathcal{G}$ be an $r$-discrete measured groupoid on $(X, \mu)$. Then the isotropy groupoid
\[\mathcal{G}_{0}=\{ \gamma \in \mathcal{G}: r(\gamma)=s(\gamma)\}\]
is normal in $\mathcal{G}$.
\end{lem}

\begin{lem}\label{lem-map-ext}
Let $\mathcal{G}$ be an $r$-discrete measured  groupoid on $(X, \mu)$ and let $A$ be a Borel subset of $X$. Then we can find a Borel map $f\colon \mathcal{G}A\rightarrow \mathcal{G}$ such that
\begin{enumerate}
\item[(i)] $s(f(x))=x$ and $r(f(x))\in A$ for a.e.\ $x\in \mathcal{G}A$;
\item[(ii)] $f(x)=e_{x}\in \mathcal{G}_{x}^{x}$ for a.e.\ $x\in A$,
where $\mathcal{G}_{x}^{x}=\{ \gamma \in \mathcal{G}: r(\gamma)=s(\gamma)=x\}$ and $e_{x}$ is the identity element of the isotropy group $\mathcal{G}_{x}^{x}$. 
\end{enumerate}
\end{lem}

\begin{lem}
Let $\mathcal{G}$ be an $r$-discrete measured  groupoid on $(X, \mu)$ and let $\mathcal{S}$ be a normal subgroupoid of $\mathcal{G}$. If $A$ is a Borel subset of $X$ with positive measure, then $(\mathcal{S})_{A}$ is normal in $(\mathcal{G})_{A}$.
\end{lem}

\begin{proof}
Let $\{ \phi_{n}\}$ be a family of normal choice functions for the pair $(\mathcal{G}, \mathcal{S})$. We write $B=\mathcal{S}A$ and
\[D_{n}=\{ x\in A\cap {\rm dom}(\phi_{n}): r(\phi_{n}(x))\in B\}.\]
Define a Borel map $\phi_{n}'\colon D_{n}\rightarrow (\mathcal{G})_{A}$ by $\phi_{n}'(x)=f(r(\phi_{n}(x)))\phi_{n}(x)$ for $x\in D_{n}$, where $f\colon B\rightarrow \mathcal{S}$ is a Borel map given by Lemma \ref{lem-map-ext} such that
\begin{itemize}
\item $s(f(x))=x$ and $r(f(x))\in A$ for a.e.\ $x\in B=\mathcal{S}A$;
\item $f(x)=e_{x}\in \mathcal{S}_{x}^{x}$ for a.e.\ $x\in A$.
\end{itemize} 

We show that $\{ \phi_{n}'\}$ is a family of normal choice functions for $((\mathcal{G})_{A}, (\mathcal{S})_{A})$. Since $\phi_{n}'$ is the composition of $\phi_{n}$ and $f$, we see that $\phi_{n}'\in {\rm End}_{\mathcal{S}}(\mathcal{G})$. Let $\gamma \in (\mathcal{G})_{A}$. Then there exists $\phi_{n}$ such that $r(\gamma)\in {\rm dom}(\phi_{n})$ and $\phi_{n}(r(\gamma))\gamma \in \mathcal{S}$. Note that $r(\gamma)\in A\cap {\rm dom}(\phi_{n})$ and $r(\phi_{n}(r(\gamma)))=r(\phi_{n}(r(\gamma))\gamma)\in \mathcal{S}A$. Therefore, $r(\gamma)\in D_{n}$ and 
\[\phi_{n}'(r(\gamma))\gamma=f(r(\phi_{n}(r(\gamma))))\phi_{n}(r(\gamma))\gamma \in (\mathcal{S})_{A},\]
which completes the proof.
\end{proof}

\begin{lem}\label{lem-normal-dir-prod}
Let $G$ be a discrete group generated by two subgroups $G_{1}$, $G_{2}$ so that $G_{1}$ is normal in $G$ and assume that we have a non-singular action of $G$ on a standard finite measure space $(X, \mu)$. We denote by $\mathcal{G}$, $\mathcal{G}_{1}$ and $\mathcal{G}_{2}$ the groupoids arising from the actions of $G$, $G_{1}$ and $G_{2}$, respectively. Let $A\subset X$ be a Borel subset with positive measure. Then $(\mathcal{G}_{1})_{A}$ is normal in the subgroupoid $\mathcal{H}=(\mathcal{G}_{1})_{A}\vee (\mathcal{G}_{2})_{A}$ of $(\mathcal{G})_{A}$ generated by the two subgroupoids $(\mathcal{G}_{1})_{A}$ and $(\mathcal{G}_{2})_{A}$.
\end{lem}

\begin{proof}
For each $i=1, 2$ and $g\in G_{i}$, define $A_{g}=A\cap g^{-1}A$ and $\psi_{g}\colon A_{g}\rightarrow (\mathcal{G})_{A}$ by $\psi_{g}(x)=(g, x)$ for $x\in A_{g}$. It is easy to check that $\psi_{g}\in {\rm End}_{(\mathcal{G}_{1})_{A}}(\mathcal{H})$. For each word $\omega$ of elements in $G_{1}$ and $G_{2}$, we can naturally define the composition $\psi_{\omega}\in {\rm End}_{(\mathcal{G}_{1})_{A}}(\mathcal{H})$ of $\psi_{g}$. It is clear that $\{ \psi_{\omega}\}_{\omega}$ forms a family of normal choice functions for $(\mathcal{H}, (\mathcal{G}_{1})_{A})$. 
\end{proof}


\section{Actions of some discrete groups}\label{sec-gr-gen}

\subsection{Actions of the mapping class group}

In Sections \ref{sec-cha-red} and \ref{sec-eq-ME}, we consider mainly the groupoid generated by an action of the mapping class group and its subgroupoids. In this subsection, we collect some results about them. Most of the following results can be shown in the same way as in \cite{kida}, where we assume that the action is essentially free.

\begin{defn}
An $r$-discrete measured groupoid $\mathcal{G}$ on $(X, \mu)$ is said to be of infinite type if there exists a Borel partition $X=A_{1}\sqcup A_{2}$ such that 
\begin{enumerate}
\item[(i)] for a.e.\ $x\in A_{1}$, the isotropy group $\mathcal{G}_{x}^{x}$ is infinite;
\item[(ii)] the associated principal groupoid of $(\mathcal{G})_{A_{2}}$ defined by
\[\{ (r(\gamma), s(\gamma))\in A_{2}\times A_{2}: \gamma \in (\mathcal{G})_{A_{2}}\}\]
is recurrent. 
\end{enumerate}
\end{defn}
Note that for any $n\in \mathbb{N}\cup \{ \infty \}$, the subset
\[X_{n}=\{ x\in X: |\mathcal{G}_{x}^{x}|=n\}\]
is measurable and satisfies $\mathcal{G}X_{n}=X_{n}$.

Let $\mathcal{G}$ be an $r$-discrete measured groupoid on $(X, \mu)$ and let $\rho \colon \mathcal{G}\rightarrow G$ be a groupoid homomorphism into a standard Borel group $G$. Let $S$ be a Borel $G$-space. Recall that a Borel map $\varphi \colon A\rightarrow S$ from a Borel subset $A\subset X$ is said to be $\rho$-invariant (or invariant for simplicity) for $\mathcal{G}$ if $\rho(\gamma)\varphi(s(\gamma))=\varphi(r(\gamma))$ for a.e.\ $\gamma \in (\mathcal{G})_{A}$.

\begin{lem}\label{lem-ext-inv}
In the above situation, let us define a Borel map $\varphi'\colon \mathcal{G}A\rightarrow S$ by $\varphi'(x)=\rho(f(x)^{-1})\varphi(r(f(x)))$ for $x\in \mathcal{G}A$, where $f\colon \mathcal{G}A\rightarrow \mathcal{G}$ is a Borel map constructed in Lemma \ref{lem-map-ext}. Then $\varphi'$ is also $\rho$-invariant for $\mathcal{G}$.
\end{lem}

\begin{proof}
Let $\gamma \in (\mathcal{G})_{\mathcal{G}A}$ and put $y=r(\gamma)$, $x=s(\gamma)\in \mathcal{G}A$. Then
\begin{align*}
\rho(\gamma)\varphi'(x)&=\rho(f(y)^{-1})\rho(f(y)\gamma f(x)^{-1})\varphi(r(f(x)))\\
&=\rho(f(y)^{-1})\varphi(r(f(y)))=\varphi'(y)
\end{align*}
since $f(y)\gamma f(x)^{-1}\in (\mathcal{G})_{A}$ and $r(f(y)\gamma f(x)^{-1})=r(f(y))$.
\end{proof}

\begin{assumption}\label{assumption-star}
We call the following assumption $(\star)$: let $\Gamma$ be a subgroup of $\Gamma(M; m)$, where $M$ is a surface with $\kappa(M)>0$ and $m\geq 3$ is an integer. Let $(X, \mu)$ be a standard finite measure space. We assume that there exists a non-singular action of $\Gamma$ on $(X, \mu)$ generating the groupoid 
\[\mathcal{G}=\mathcal{G}_{\Gamma}=\{ (\gamma, x)\in \Gamma \times X: \gamma \in \Gamma,\ x\in X\}.\]
and the induced cocycle
\[\rho \colon \mathcal{G}\rightarrow \Gamma,\ \ (\gamma, x)\mapsto \gamma\]
for $\gamma \in \Gamma$ and a.e.\ $x\in X$. 
\end{assumption}

Under the above assumption, we often use the following notation: 
\begin{itemize}
\item For a subgroup $\Gamma'$ of $\Gamma$, let $\mathcal{G}_{\Gamma'}$ denote the subgroupoid of $\mathcal{G}$ generated by the action of $\Gamma'$.  

\item For $\sigma \in S(M)$, we denote by $D_{\sigma}$ the intersection of $\Gamma$ and the subgroup generated by Dehn twists about all curves in $\sigma$. We write $\mathcal{G}_{\sigma}$ instead of $\mathcal{G}_{D_{\sigma}}$ for simplicity. If $\sigma$ consists of one element $\alpha \in V(C)$, then we write $D_{\alpha}$ (resp. $\mathcal{G}_{\alpha}$) instead of $D_{\sigma}$ (resp. $\mathcal{G}_{\sigma}$).
\end{itemize}

As in \cite{kida}, we can consider two types of subgroupoids of infinite type, following the classification of subgroups of $\Gamma(M)$ mentioned in Subsection \ref{subsec-mcg}. We denote by $M(\mathcal{PMF})$ the space of all probability measures on $\mathcal{PMF}$.

\begin{thm}[\ci{Theorem 4.41}{kida}]\label{thm-alternative}
Under the assumption $(\star)$, let $Y\subset X$ be a Borel subset with positive measure and let $\mathcal{S}$ be a subgroupoid of $(\mathcal{G})_{Y}$ of infinite type. If we have an invariant Borel map $\varphi \colon Y\rightarrow M(\mathcal{PMF})$ for $\mathcal{S}$, then there exists a Borel partition $Y=Y_{1}\sqcup Y_{2}$ satisfying the following:
\begin{enumerate}
\item[(i)] $\varphi(x)(\mathcal{MIN})=1$ for a.e.\ $x\in Y_{1}$;
\item[(ii)] $\varphi(x)(\mathcal{PMF}\setminus \mathcal{MIN})=1$ for a.e.\ $x\in Y_{2}$.
\end{enumerate}
\end{thm}

In the above theorem, remark that both $Y_{1}$ and $Y_{2}$ are invariant for $\mathcal{S}$ and if $Y'$ is a Borel subset of $Y$ with positive measure and $\psi \colon Y'\rightarrow M(\mathcal{PMF})$ is another invariant Borel map for $\mathcal{S}$, then $\psi$ satisfies
\begin{enumerate}
\item[(i)] $\psi(x)(\mathcal{MIN})=1$ for a.e.\ $x\in Y_{1}\cap Y'$;
\item[(ii)] $\psi(x)(\mathcal{PMF}\setminus \mathcal{MIN})=1$ for a.e.\ $x\in Y_{2}\cap Y'$,
\end{enumerate}
where $Y_{1}$, $Y_{2}$ are Borel subsets as in the above theorem. Hence, it is natural to give the following definition:

\begin{defn}\label{defn-ia-red}
Under the assumption $(\star)$, let $Y\subset X$ be a Borel subset with positive measure and let $\mathcal{S}$ be a subgroupoid of $(\mathcal{G})_{Y}$ of infinite type.
\begin{enumerate}  
\item[(i)] If we have an invariant Borel map $\varphi \colon Y\rightarrow M(\mathcal{PMF})$ for $\mathcal{S}$ such that $\varphi(x)(\mathcal{MIN})=1$ for a.e.\ $x\in Y$, then $\mathcal{S}$ is said to be irreducible and amenable (or IA in short). 
\item[(ii)] If we have an invariant Borel map $\varphi \colon Y\rightarrow M(\mathcal{PMF})$ for $\mathcal{S}$ such that $\varphi(x)(\mathcal{PMF}\setminus \mathcal{MIN})=1$ for a.e.\ $x\in Y$, then $\mathcal{S}$ is said to be reducible.
\end{enumerate}
\end{defn}

We have explained in Section \ref{sec-int} that the key ingredient of the proof of Theorem \ref{thm-main} is to construct an essentially unique almost $\Gamma(M)\times \Gamma(M)$-equivariant Borel map $\Phi \colon \Sigma \rightarrow {\rm Aut}(C)$ for a self ME coupling $(\Sigma, m)$ of $\Gamma(M)$. We give a rough outline of the construction of the map $\Phi$. In \cite{kida}, we developed the theory of recurrent subrelations of an equivalence relation arising from a standard action of the mapping class group. Thanks to it, we can divide such subrelations into two types, irreducible and amenable (IA) ones and reducible ones as in Theorem \ref{thm-alternative} and Definition \ref{defn-ia-red}. The notion of normal subrelations also played an important role in the classification theorem of \cite{kida}. In what follows in this section, we generalize various central results in \cite{kida} about the above subrelations to the case where the action of the mapping class group is not necessarily essentially free. In this case, although we need to consider $r$-discrete measured groupoids arising from group actions, the proof can be proceeded along the same line.

In Section \ref{sec-cha-red}, using various results in this section, we characterize a reducible subgroupoid in terms of amenability, non-amenability and normality (see Propositions \ref{prop-cha-red-ame}, \ref{prop-cha-red-non}). Note that these three properties are preserved under an isomorphism between two groupoids, and that as mentioned in Subsection \ref{subsec-ME oe}, considering a self ME coupling of $\Gamma(M)$ is almost equivalent to considering an isomorphism $f$ between two groupoids $\mathcal{G}_{1}$, $\mathcal{G}_{2}$ arising from measure-preserving actions of $\Gamma(M)$. Thanks to the characterization of a reducible subgroupoid, we see that the image of a reducible subgroupoid of $\mathcal{G}_{1}$ via $f$ is also reducible. Moreover,  maximal reducible subgroupoids are mapped to maximal ones by $f$ (see Corollary \ref{cor-max-red}, Lemma \ref{lem-max-map} and Corollary \ref{cor-pre-max-map}).

Let $\mathcal{G}$ be the groupoid associated with a measure-preserving action of $\Gamma(M)$ on a standard finite measure space $(X, \mu)$. As a next stage, in Section \ref{sec-eq-ME}, we study an amenable normal subgroupoid $\mathcal{S}$ of infinite type of a maximal reducible subgroupoid of $\mathcal{G}$. We can show that $\mathcal{S}$ is contained in the groupoid generated by the action of a Dehn twist about some simple closed curve on $M$ up to a countable Borel partition of $X$ (see Lemma \ref{lem-dehn-max-normal}). Conversely, the subgroupoid generated by the action of a Dehn twist about a simple closed curve is normal in some maximal reducible subgroupoid of $\mathcal{G}$. It follows that roughly speaking, the subgroupoid generated by a Dehn twist about a simple closed curve can be characterized in terms of amenability, non-amenability and normality, and in the situation of the previous paragraph, we see that such subgroupoids are preserved by $f$. This means that $f$ induces a bijection of the set $V(C)$ of all isotopy classes of simple closed curves on $M$, which is shown to be an automorphism of the curve complex. Translating this fact into structural information on a self ME coupling $(\Sigma, m)$ of $\Gamma(M)$, we can construct an almost $\Gamma(M)\times \Gamma(M)$-equivariant Borel map $\Phi$ from $\Sigma$ into ${\rm Aut}(C)$ as mentioned in Section \ref{sec-int}.

\begin{rem}
When $\kappa(M)>0$ and $M\neq M_{1, 2}, M_{2, 0}$, Ivanov \cite[Section 8.5]{ivanov2} showed that any isomorphism from a finite index subgroup of $\Gamma(M)^{\diamond}$ into a finite index subgroup of $\Gamma(M)^{\diamond}$ is the restriction of a unique inner automorphism of $\Gamma(M)^{\diamond}$. A key ingredient of his proof is to show that such an isomorphism $f$ maps sufficiently high powers of Dehn twists into powers of Dehn twists by characterizing a (power of) Dehn twist algebraically (see \cite[Theorem 7.5.B]{ivanov2}). It follows that $f$ yields a bijection on the set $V(C)$, which is, in fact, an automorphism of the curve complex. This automorphism comes from an element $g$ in $\Gamma(M)^{\diamond}$ by Ivanov's theorem (see Theorem \ref{thm-cc-auto}). After an easy computation shown in the proof of Theorem 8.5.A in \cite{ivanov2}, we can prove that $f$ is the restriction of the inner automorphism of $\Gamma(M)$ by conjugation with $g$. Our construction of the map $\Phi$ mentioned above relies heavily on this idea due to Ivanov. 
\end{rem}

\begin{rem}
Note that if we only want to prove Theorem \ref{thm-main}, then it is not necessary to generalize the results in \cite{kida} to the case where the standard action of the mapping class group is not necessarily essentially free because in general, when two discrete groups $\Lambda_{1}$, $\Lambda_{2}$ are measure equivalent, there exists a ME coupling of $\Lambda_{1}$ and $\Lambda_{2}$ such that the $\Lambda_{1}\times \Lambda_{2}$-action on it is essentially free, which induces weakly orbit equivalence between standard actions of $\Lambda_{1}$ and $\Lambda_{2}$. However, in Theorem \ref{thm-int-lattice-emb}, we need to consider a ME coupling of $\Gamma(M)$ such that the $\Gamma(M)\times \Gamma(M)$-action is not necessarily essentially free. Moreover, thanks to the generalization, we can obtain some information about stabilizers of measure-preserving actions of $\Gamma(M)$ (see Corollaries \ref{cor-mcg-stab} and \ref{cor-hyp-stab}). 
\end{rem}

In the following theorems, we collect basic properties of the two types of subgroupoids of a groupoid generated by an action of the mapping class group, which are called IA and reducible subgroupoids, respectively. First, we consider IA subgroupoids.

Since the curve complex $C$ is a hyperbolic metric space \cite{masur-minsky1}, we can construct the boundary $\partial C$ at infinity, which is non-empty (see \cite{ham}, \cite{kla} or \cite[Section 3.2]{kida}). Let $\partial_{2}C$ be the quotient space of $\partial C\times \partial C$ by the coordinate exchanging action of the symmetric group of two letters and let $M(\partial C)$ be the space of all probability measures on $\partial C$, which has the Borel structure introduced in the comment before Proposition 4.30 in \cite{kida}. Each element in $\partial_{2}C$ can naturally be viewed as an atomic measure in $M(\partial C)$ so that each atom has measure $1$ or $1/2$. Then $\partial_{2}C$ is a Borel subset of $M(\partial C)$.

Under the assumption $(\star)$, let $Y\subset X$ be a Borel subset with positive measure and let $\mathcal{S}$ be a subgroupoid of $(\mathcal{G})_{Y}$ of infinite type. Note that if $\mathcal{S}$ is IA, then we can construct an invariant Borel map $Y\rightarrow M(\partial C)$ for $\mathcal{S}$ by using the $\Gamma(M)$-equivariant map $\mathcal{MIN}\rightarrow \partial C$, which is constructed in \cite{kla} (see also \cite{ham}).

\begin{prop}[\ci{Proposition 4.32 (ii), Corollary 4.43}{kida}]\label{prop-bou-IA}
Under the assumption $(\star)$, let $Y\subset X$ be a Borel subset with positive measure and let $\mathcal{S}$ be a subgroupoid of $(\mathcal{G})_{Y}$ of infinite type. Suppose that $\mathcal{S}$ admits a $\rho$-invariant Borel map $\varphi \colon Y\rightarrow M(\partial C)$. Then the cardinality of ${\rm supp}(\varphi(x))$ is at most $2$ for a.e.\ $x\in Y$ and $\mathcal{S}$ is IA. 
\end{prop}

\begin{thm}[\ci{Section 4.4.1, Lemma 4.58}{kida}]\label{thm-ia}
Under the assumption $(\star)$, let $Y\subset X$ be a Borel subset with positive measure and let $\mathcal{S}$ be a subgroupoid of $(\mathcal{G})_{Y}$ of infinite type. Suppose that $\mathcal{S}$ is IA. Then
\begin{enumerate}
\item[(i)] there exists an essentially unique invariant Borel map $\varphi_{0}\colon Y\rightarrow \partial_{2}C$ for $\mathcal{S}$ satisfying the following: if $Y'$ is a Borel subset of $Y$ with positive measure and $\varphi \colon Y'\rightarrow M(\partial C)$ is an invariant Borel map for $\mathcal{S}$, then
\[{\rm supp}(\varphi(x))\subset {\rm supp}(\varphi_{0}(x))\]
for a.e.\ $x\in Y'$, where ${\rm supp}(\nu)$ denotes the support of a measure $\nu$.
\item[(ii)] if $\mathcal{T}$ is a subgroupoid of $(\mathcal{G})_{Y}$ with $\mathcal{S}\vartriangleleft \mathcal{T}$, then $\varphi_{0}$ is invariant also for $\mathcal{T}$.
\item[(iii)] the groupoid $\mathcal{S}$ is amenable. 
\end{enumerate}
\end{thm}

If $X$ is a point and $\mathcal{G}$ is isomorphic to $\Gamma$, then the above facts follow from the classification of subgroups of $\Gamma(M)$ described in Subsection \ref{subsec-mcg}. In this case, using properties of pseudo-Anosov elements, we can prove that $\mathcal{S}$ is virtually cyclic, which implies Theorem \ref{thm-ia} (iii). To prove Theorem \ref{thm-ia} (iii) in a general case, we need to use the amenability in a measurable sense of the action of $\Gamma(M)$ on $\partial C$ (or $\partial_{2}C$).

\begin{rem}\label{rem-div}
Under the assumption $(\star)$, let $Y\subset X$ be a Borel subset with positive measure and let $\mathcal{S}$ be a subgroupoid of $(\mathcal{G})_{Y}$ of infinite type. It follows from Theorem \ref{thm-alternative} that there exists an essentially unique Borel partition $Y=Y_{1}\sqcup Y_{2}\sqcup Y_{3}$ satisfying the following:
\begin{itemize}
\item if $Y_{1}$ has positive measure, then $(\mathcal{S})_{Y_{1}}$ is IA;
\item if $Y_{2}$ has positive measure, then $(\mathcal{S})_{Y_{2}}$ is reducible;
\item if $Y_{3}$ has positive measure, then $(\mathcal{S})_{Y_{3}'}$ admits no invariant Borel maps $Y_{3}'\rightarrow M(\mathcal{PMF})$ for any Borel subset $Y_{3}'$ of $Y_{3}$ with positive measure. 
\end{itemize}
If $\mathcal{S}$ is amenable and any restriction of $\mathcal{S}$ to a Borel subset of $Y$ with positive measure is not reducible, then $\mathcal{S}$ is IA, and the converse also holds by Theorem \ref{thm-ia}. 
\end{rem}

Next, we recall basic properties of reducible subgroupoids. We can define the canonical reduction system for a reducible subgroupoid as in the case of a reducible subgroup.

\begin{defn}
Under the assumption $(\star)$, let $Y\subset X$ be a Borel subset with positive measure and let $\mathcal{S}$ be a subgroupoid of $(\mathcal{G})_{Y}$ of infinite type. Let $A$ be a Borel subset of $Y$ with positive measure and let $\alpha \in V(C)$. 
\begin{enumerate}
\item[(i)] We say that the pair $(\alpha, A)$ is $\rho$-invariant for $\mathcal{S}$ if there exists a countable Borel partition $A=\bigsqcup A_{n}$ of $A$ such that the constant map $A_{n}\rightarrow \{ \alpha \}$ is invariant for $\mathcal{S}$ for each $n$. 
\item[(ii)] Suppose that $(\alpha, A)$ is $\rho$-invariant for $\mathcal{S}$. The pair $(\alpha, A)$ is said to be purely $\rho$-invariant if $(\beta, B)$ is not $\rho$-invariant for $\mathcal{S}$ for any Borel subset $B$ of $A$ with positive measure and any $\beta \in V(C)$ with $i(\alpha, \beta)\neq 0$. (In \cite{kida}, we call this pair an essential $\rho$-invariant one for $\mathcal{S}$.)
\end{enumerate} 
\end{defn}

Since we consider a countable Borel partition in the definition of a $\rho$-invariant pair, we can easily show that if there are $\alpha \in V(C)$ and a Borel subset $A_{n}\subset Y$ for $n\in \mathbb{N}$ with $(\alpha, A_{n})$ $\rho$-invariant for $\mathcal{S}$ (or purely $\rho$-invariant, respectively), then the pair $(\alpha, \bigcup A_{n})$ is also $\rho$-invariant for $\mathcal{S}$ (or purely $\rho$-invariant, respectively) by \cite[Lemma 4.47]{kida}. It follows that for each $\alpha \in V(C)$, we can find an essentially maximal Borel subset $A_{\alpha}\subset Y$ such that $(\alpha, A_{\alpha})$ is $\rho$-invariant for $\mathcal{S}$ (or purely $\rho$-invariant, respectively) unless there exist no $\rho$-invariant pairs for $\mathcal{S}$ (or purely $\rho$-invariant, respectively).

\begin{thm}[\ci{Section 4.5, Lemma 4.60}{kida}]\label{thm-reducible}     
Under the assumption $(\star)$, let $Y\subset X$ be a Borel subset with positive measure and let $\mathcal{S}$ be a subgroupoid of $(\mathcal{G})_{Y}$ of infinite type. Suppose that $\mathcal{S}$ is reducible. Then
\begin{enumerate}
\item[(i)] there exists a purely $\rho$-invariant pair for $\mathcal{S}$.
\item[(ii)] we can define an essentially unique invariant Borel map $\varphi \colon Y\rightarrow S(M)$ for $\mathcal{S}$ so that 
\begin{enumerate}
\item[(a)] if $\sigma \in S(M)$ satisfies $\mu(\varphi^{-1}(\sigma))>0$ and $\alpha \in \sigma$, then $(\alpha, \varphi^{-1}(\sigma))$ is a purely $\rho$-invariant pair for $\mathcal{S}$;
\item[(b)] if $(\alpha, A)$ is a purely $\rho$-invariant pair for $\mathcal{S}$, then 
\[\mu(A\setminus \varphi^{-1}(\{ \sigma \in S(M): \alpha \in \sigma \}))=0.\]
\end{enumerate} 
\item[(iii)] if $\mathcal{T}$ is a subgroupoid of $(\mathcal{G})_{Y}$ with $\mathcal{S}\vartriangleleft \mathcal{T}$, then $\varphi$ is invariant also for $\mathcal{T}$.
\end{enumerate}
\end{thm}

Note that we can construct $\varphi$ in Theorem \ref{thm-reducible} (ii) from (i) by the observation right before Theorem \ref{thm-reducible}. We call $\varphi$ in the above theorem the canonical reduction system (CRS) for $\mathcal{S}$. It is easy to see that if $A$ is a Borel subset of $Y$ with positive measure, then the CRS for $(\mathcal{S})_{A}$ is the restriction of $\varphi$ to $A$ (see \cite[Lemma 4.53 (iii)]{kida}). If $X$ is a point and $\mathcal{G}$ is isomorphic to $\Gamma$, then the above definition of the CRS for $\mathcal{S}$ coincides with the one mentioned in Subsection \ref{subsec-mcg}. 

In the following two lemmas, we study the CRS's for certain reducible subgroupoids arising from the actions of reducible subgroups.

\begin{lem}\label{lem-red-crs}
Under the assumption $(\star)$, suppose that the $\Gamma$-action on $(X, \mu)$ is measure-preserving. Let $G$ be an infinite reducible subgroup of $\Gamma$ and let $\sigma \in S(M)$ be the CRS for $G$. Then $\mathcal{G}_{G}$ is reducible and its CRS $\varphi \colon X\rightarrow S(M)$ is constant with value $\sigma$. 
\end{lem}

\begin{proof}
It is clear that $\mathcal{G}_{G}$ is reducible and for any $\alpha \in \sigma$, the pair $(\alpha, X)$ is $\rho$-invariant for $\mathcal{G}_{G}$. Assume that there exists $\alpha \in \sigma$ such that the pair $(\alpha, X)$ is not purely $\rho$-invariant for $\mathcal{G}_{G}$. Then we would have a Borel subset $B$ of $X$ with positive measure and $\beta \in V(C)$ with $i(\alpha, \beta)\neq 0$ such that $(\beta, B)$ is a $\rho$-invariant pair for $\mathcal{G}_{G}$. It follows that there exists a Borel subset $B'$ of $B$ with positive measure such that $\rho(\gamma)\beta =\beta$ for a.e.\ $\gamma \in (\mathcal{G}_{G})_{B'}$. We can find $g\in G$ of infinite order with $g\beta \neq \beta$ since $\alpha \in \sigma$. Since $g$ has infinite order and the $\Gamma$-action on $(X, \mu)$ preserves the finite positive measure $\mu$, the subgroupoid $(\mathcal{G}_{\langle g\rangle})_{B'}$ is of infinite type, where $\langle g\rangle$ denotes the cyclic subgroup generated by $g$. We can find a Borel subset $B_{1}'\subset B'$ with positive measure and $n\in \mathbb{Z}\setminus \{ 0\}$ such that $(g^{n}, x)\in (\mathcal{G}_{\langle g\rangle })_{B'}$ for any $x\in B_{1}'$. Thus, $g^{n}\beta =\rho(g^{n}, x)\beta =\beta$ holds for any $x\in B_{1}'$. Since $G$ is a subgroup of $\Gamma(M; m)$, it follows from Theorem \ref{thm-pure} that $g^{k}\beta =\beta$ for any $k\in \mathbb{Z}\setminus \{ 0\}$, which is a contradiction. Thus, $(\alpha, X)$ is a purely $\rho$-invariant pair for $\mathcal{G}_{G}$ and we see that $\sigma$ is contained in $\varphi(x)$ for a.e.\  $x\in X$.

Next, assume that we have a Borel subset $A$ of $X$ with positive measure and $\beta \in \varphi(x)\setminus \sigma$ for any $x\in A$. For each $g\in G$, there is a Borel subset $A_{1}\subset A$ with positive measure and $n\in \mathbb{Z}\setminus \{ 0\}$ such that $(g^{n}, x)\in (\mathcal{G}_{G})_{A_{1}}$ and the equation $g^{n}\beta =\rho(g^{n}, x)\beta =\beta$ holds for any $x\in A_{1}$ because $(\beta, A)$ is a $\rho$-invariant pair for $\mathcal{G}_{G}$. Thus, $g\beta =\beta$ for any $g\in G$ by Theorem \ref{thm-pure}. It follows from $\beta \not\in \sigma$ that there exists $\gamma \in V(C)$ such that $i(\beta, \gamma)\neq 0$ and $h\gamma =\gamma$ for any $h\in G$. Thus, $(\gamma, X)$ is a $\rho$-invariant pair for $\mathcal{G}_{G}$, which contradicts the assumption that $\beta \in \varphi(x)$ for any $x\in A$, that is, the pair $(\beta, A)$ is pure.    
\end{proof}

\begin{lem}\label{lem-crs-dehn}
Under the assumption $(\star)$, suppose that the $\Gamma$-action on $(X, \mu)$ is measure-preserving. Let $\alpha \in V(C)$ and assume that the subgroup $D_{\alpha}$ is infinite. Let $Y$ be a Borel subset of $X$ with positive measure and let $\mathcal{S}$ be a subgroupoid of $(\mathcal{G}_{\alpha})_{Y}$ of infinite type. Then $\mathcal{S}$ is reducible and its CRS for $\mathcal{S}$ is constant with value $\{ \alpha \}$.
\end{lem}

\begin{proof}
It is clear that $\mathcal{S}$ is reducible and the pair $(\alpha, Y)$ is $\rho$-invariant for $\mathcal{S}$. Let $A$ be a Borel subset of $Y$ with positive measure and $\beta \in V(C)$ with $i(\alpha, \beta)\neq 0$. Assume that the pair $(\beta, A)$ is $\rho$-invariant for $\mathcal{S}$. Then there exists a Borel subset $B$ of $A$ with positive measure such that $\rho(\gamma)\beta =\beta$ for a.e.\ $\gamma \in (\mathcal{S})_{B}$. Since $\mathcal{S}$ is a subgroupoid of $(\mathcal{G}_{\alpha})_{Y}$ of infinite type, there exist infinitely many $n\in \mathbb{Z}$ and a Borel subset $B_{n}$ of $B$ with positive measure such that $t^{n}\in \Gamma$ and $(t^{n}, x)\in (\mathcal{S})_{B}$ for any $x\in B_{n}$, where $t\in \Gamma(M)$ denotes a Dehn twist about $\alpha$. Hence, $t^{n}\beta =\rho(t^{n}, x)\beta =\beta$ for a.e.\ $x\in B_{n}$. In particular, $t^{n}\beta =\beta$ for infinitely many $n\in \mathbb{Z}$. It follows from \cite[Lemma 4.2]{ivanov1} (or \cite[Corollary 4.26]{kida}) that $i(\alpha, \beta)=0$, which is a contradiction. Thus, the pair $(\alpha, Y)$ is a pure $\rho$-invariant one for $\mathcal{S}$.   

If $\gamma \in V(C)$ satisfies $i(\alpha, \gamma)=0$ and $\alpha \neq \gamma$, then there exists $\delta \in V(C)$ such that $i(\alpha, \delta)=0$ and $i(\gamma, \delta)\neq 0$. Since the pair $(\delta, Y)$ is $\rho$-invariant for $\mathcal{S}$, the pair $(\gamma, A')$ can not be a pure $\rho$-invariant one for $\mathcal{S}$ for any Borel subset $A'$ of $Y$.
\end{proof}

The following proposition can also be proved along the same line as in \cite{kida}:

\begin{prop}[\ci{Proposition 4.61}{kida}]\label{prop-not-ia-red}
Under the assumption $(\star)$, suppose that the $\Gamma$-action on $(X, \mu)$ is measure-preserving and that $\Gamma$ is sufficiently large. Then $(\mathcal{G})_{Y}$ is neither IA nor reducible for any Borel subset $Y\subset X$ with positive measure. 
\end{prop}

As an application of the above generalization of the results in \cite{kida}, we obtain some information on stabilizers for a measure-preserving action of the mapping class group on a standard finite measure space.

\begin{cor}\label{cor-mcg-stab}
Under the assumption $(\star)$, suppose that the $\Gamma$-action on $(X, \mu)$ is measure-preserving and that $\Gamma$ is sufficiently large. Then the isotropy group 
\[\mathcal{G}_{x}^{x}=\{ \gamma \in \mathcal{G}: r(\gamma)=s(\gamma)=x\}\]
is either trivial or sufficiently large for a.e.\ $x\in X$.
\end{cor}

We can show this corollary by using Lemma \ref{lem-iso-normal}, Theorem \ref{thm-ia} (ii), Theorem \ref{thm-reducible} (iii) and Proposition \ref{prop-not-ia-red}. Note that $\Gamma$ is torsion-free and that
\begin{itemize}
\item for each pseudo-Anosov element $g\in \Gamma$, the subset of $X$ consisting of $x\in X$ such that $\mathcal{G}_{x}^{x}$ is IA and fixes the pair $\{ F_{\pm}(g)\}$ of pseudo-Anosov foliations is measurable; 
\item for each $\sigma \in S(M)$, the subset of $X$ consisting of $x\in X$ such that $\mathcal{G}_{x}^{x}$ is reducible and its CRS is $\sigma$ is measurable.
\end{itemize}
It follows from these remarks that for a Borel subset $Y$ of $X$ with positive measure, both subsets
\[Y_{1}=\{ x\in Y: \mathcal{G}_{x}^{x} {\rm \ is\ IA}\},\ \ Y_{2}=\{ x\in Y: \mathcal{G}_{x}^{x} {\rm \ is\ reducible}\}\]
are measurable, and $(\mathcal{G}_{0})_{Y_{1}}$ is IA and $(\mathcal{G}_{0})_{Y_{2}}$ is reducible, where 
\[\mathcal{G}_{0}=\{ \gamma \in \mathcal{G}: r(\gamma)=s(\gamma)\}\]
is the isotropy groupoid of $\mathcal{G}$.

In order to analyze reducible subgroupoids furthermore, in Theorems \ref{thm-alternative-trivial-irreducible}, \ref{thm-alternative-irreducible-amenable-non-amenable}, we consider components of the surface obtained by cutting along the CRS for $\mathcal{S}$. There are three types of components as in the case of subgroups of $\Gamma(M; m)$ mentioned in the comment right before Lemma \ref{lem-red-crs-group}. These theorems will be used when we characterize reducible subgroupoids in terms of amenability, non-amenability and normality (see Section \ref{sec-cha-red}).

If $\Gamma$ is an infinite reducible subgroup of $\Gamma(M; m)$ with an integer $m\geq 3$ and $\sigma \in S(M)$ is the CRS for $\Gamma$, then we can classify each component $Q$ of $M_{\sigma}$ in terms of properties of the quotients $p_{Q}(\Gamma)$ (i.e., triviality and amenability), where $p_{Q}\colon \Gamma \rightarrow \Gamma(Q)$ is the natural homomorphism (see the comment right before Lemma \ref{lem-red-crs-group}). On the other hand, when we consider a reducible subgroupoid, we cannot construct such a quotient. However, fortunately, the properties of the quotient $p_{Q}(\Gamma)$ used in the classification of $Q$ can be characterized in terms of fixed points for the action of $p_{Q}(\Gamma)$ on the space $M(\mathcal{PMF}(Q))$ of all probability measures on $\mathcal{PMF}(Q)$ as follows:
\begin{enumerate}
\item[(a)] $Q$ is T for $\Gamma$ if and only if either $Q$ is a pair of pants or $p_{Q}(g)\alpha =\alpha$ for any $g\in \Gamma$ and any (or some) $\alpha \in V(C(Q))$.
\item[(b)] $Q$ is IA for $\Gamma$ if and only if the following three conditions are satisfied:
\begin{itemize}
\item $Q$ is not a pair of pants;
\item $p_{Q}(g)\alpha \neq \alpha$ for any non-trivial $g\in \Gamma$ and any (or some) $\alpha \in V(C(Q))$;
\item there exists $\mu \in M(\mathcal{PMF}(Q))$ such that $p_{Q}(g)\mu =\mu$ for any $g\in \Gamma$ and $\mu(\mathcal{MIN}(Q))=1$.
\end{itemize}
\item[(c)] $Q$ is IN for $\Gamma$ if and only if the following two conditions are satisfied:
\begin{itemize}
\item $Q$ is not a pair of pants;
\item there exist no fixed points for the action of $p_{Q}(\Gamma)$ on $M(\mathcal{PMF}(Q))$.
\end{itemize}
\end{enumerate}
By this observation, we can similarly consider three types of components of the surface obtained by cutting along the CRS for a reducible subgroupoid. Before stating the definition of the three types of components, we recall some notation.

Let $L$ be a submanifold of the surface $M$ which is a realization of some element in $S(M)$. Let $Q$ be a component of $M_{L}$, where $M_{L}$ denotes the surface obtained by cutting $M$ along $L$. Let $p_{L}\colon M_{L}\rightarrow M$ denote the canonical map. For $\delta \in V(C(M))$, we define a finite subset $r(\delta, Q)$ of $V(C(Q))$.

Let $\delta \in V(C(M))$ and represent the isotopy class $\delta$ by a circle $D$ that intersects each of the components of $L$ in the least possible number of points. Put $D_{L}=p_{L}^{-1}(D)$. The manifold $D_{L}$ consists of some intervals or it is a circle (if $D\cap L=\emptyset$). 

If either $D_{L}\cap Q=\emptyset$ or $D_{L}$ is a circle which lies in $Q$ and is peripheral for $Q$, then put $r(\delta, Q)=\emptyset$. 
If $D_{L}$ is a non-peripheral circle lying in $Q$, put $r(\delta, Q)=\{\delta \}$.

In the remaining cases, the intersection $D_{L}\cap Q$ consists of some intervals. For each such interval $I$, consider a regular neighborhood in $Q$ of the union of the interval $I$ and those components of $\partial Q$ on which the ends of $I$ lie. Let $N_{I}$ denote the regular neighborhood. Then $N_{I}$ is a disk with two holes. Let $r'(\delta, Q)$ be the set of isotopy classes of components of the manifolds $\partial N_{I}\setminus \partial Q$, where $I$ runs through the set of all components of $D_{L}\cap Q$. Define $r(\delta, Q)$ as the resulting set of discarding from $r'(\delta, Q)$ the isotopy classes of trivial or peripheral circles of $Q$. We will regard $r(\delta, Q)$ as a subset of $V(C(M))$ using the embedding $V(C(Q))\hookrightarrow V(C(M))$. It is clear that this definition depends only on $\delta$ and the isotopy class of $Q$.

Let $F\colon M\rightarrow M$ be a diffeomorphism such that $F(L)=L$ and the induced diffeomorphism $M_{L}\rightarrow M_{L}$ takes $Q$ to $Q$. If $f\in \Gamma(M)$ denotes the isotopy class of $F$, then we have the equality
\[f(r(\delta, Q))=r(f\delta, Q)\]
by definition.

\begin{lem}[\ci{Lemma 7.9}{ivanov1}]\label{lem-pants}
Let $L$ and $Q$ be the same as above and let $\delta \in V(C(M))$. If $r(\delta, Q)=\emptyset$, then one of the following three cases occurs:
\begin{enumerate}
\renewcommand{\labelenumi}{\rm(\roman{enumi})}
\item there is a circle in the class $\delta$ that does not intersect $Q$;
\item $\delta$ is the isotopy class of one of the components of $L$;
\item $Q$ is a disk with two holes.
\end{enumerate} 
\end{lem}

We denote by $D=D(M)$ the set of all isotopy classes of subsurfaces in $M$ and denote by $\mathcal{F}_{0}(D)$ the set of all finite subsets $F$ of $D$ (including the empty set) such that if $Q_{1}, Q_{2}\in F$ and $Q_{1}\neq Q_{2}$, then $Q_{1}$ and $Q_{2}$ can be realized disjointly on $M$.

\begin{thm}[\ci{Theorem 5.6}{kida}]\label{thm-alternative-trivial-irreducible}
Under the assumption $(\star)$, let $Y\subset X$ be a Borel subset with positive measure and let $\mathcal{S}$ be a subgroupoid of $(\mathcal{G})_{Y}$ of infinite type. Suppose that $\mathcal{S}$ is reducible and let $\varphi \colon Y\rightarrow S(M)$ be its CRS. Then there exist two essentially unique invariant Borel maps $\varphi_{t}, \varphi_{i}\colon Y\rightarrow \mathcal{F}_{0}(D)$ for $\mathcal{S}$ satisfying the following:
\begin{enumerate}
\item[(i)] any element in $\varphi_{t}(x)\cup \varphi_{i}(x)$ is a component of $M_{\varphi(x)}$ for a.e.\ $x\in Y$;
\item[(ii)] each component of $M_{\varphi(x)}$ belongs to $\varphi_{t}(x)\cup \varphi_{i}(x)$ and $\varphi_{t}(x)\cap \varphi_{i}(x)=\emptyset$ for a.e.\ $x\in Y$;
\item[(iii)] if $Q$ is in $F\in \mathcal{F}_{0}(D)$ with $\mu(\varphi_{t}^{-1}(F))>0$, then either $Q$ is a pair of pants or the pair $(\alpha, \varphi_{t}^{-1}(F))$ is $\rho$-invariant for $\mathcal{S}$ for any $\alpha \in V(C(Q))$;
\item[(iv)] if $Q$ is in $F\in \mathcal{F}_{0}(D)$ with $\mu(\varphi_{i}^{-1}(F))>0$, then $Q$ is not a pair of pants and $(\alpha, A)$ is not $\rho$-invariant for $\mathcal{S}$ for any $\alpha \in V(C(M))$ with $r(\alpha, Q)\neq \emptyset$ and any Borel subset $A\subset \varphi_{i}^{-1}(F)$ with positive measure.
\end{enumerate}
\end{thm}

We call $\varphi_{t}$ (resp. $\varphi_{i}$) the system of trivial (resp. irreducible) subsurfaces for $\mathcal{S}$ or in short, a T (resp. IA) system for $\mathcal{S}$. We often call an element in $\varphi_{t}(x)$ (resp. $\varphi_{i}(x)$) a trivial or T (resp. irreducible or I) subsurface for $x\in Y$. When we identify a subsurface with a component of the surface obtained by cutting along some curves, we call T and I subsurfaces also T and I components, respectively. It is easy to see that if $A$ is a Borel subset of $Y$ with positive measure, then the T, I systems for $(\mathcal{S})_{A}$ are the restrictions of $\varphi_{t}$, $\varphi_{i}$ to $A$, respectively (see \cite[Lemma 5.7]{kida}). We can show that T components have the following stronger property:

\begin{lem}[\ci{Lemma 5.4}{kida}]\label{lem-T-strong}
Under the assumption $(\star)$, let $Y\subset X$ be a Borel subset with positive measure and let $\mathcal{S}$ be a subgroupoid of $(\mathcal{G})_{Y}$ of infinite type. Suppose that $\mathcal{S}$ is reducible and let $\varphi \colon Y\rightarrow S(M)$ be its CRS. We assume the following: 
\begin{itemize}
\item $\varphi$ is constant with value $\sigma \in S(M)$ and $Q$ is a component of $M_{\sigma}$;
\item we have $\alpha \in V(C(M))$ with $r(\alpha, Q)\neq \emptyset$ and a Borel subset $A\subset Y$ with positive measure such that $(\alpha, A)$ is $\rho$-invariant for $\mathcal{S}$. 
\end{itemize}
Then there exists a countable Borel partition $A=\bigsqcup A_{n}$ such that $\rho(\gamma)\beta =\beta$ for any $\beta \in V(C(Q))$ and a.e.\ $\gamma \in (\mathcal{S})_{A_{n}}$. In particular, the pair $(\beta, A)$ is $\rho$-invariant for $\mathcal{S}$ for any curve $\beta \in V(C(Q))$.   
\end{lem}

If the cocycle $\rho \colon \mathcal{S}\rightarrow \Gamma$ is essentially valued in $\Gamma_{\sigma}=\{ g\in \Gamma :g\sigma =\sigma \}$ for some $\sigma$ and $Q$ is a component of $M_{\sigma}$, then $\rho_{Q}$ denotes the cocycle defined by the composition of $\rho$ with $p_{Q}\colon \Gamma_{\sigma}\rightarrow \Gamma(Q)$. In the next theorem, we further divide I subsurfaces into two types, IA and IN ones.

\begin{thm}[\ci{Theorem 5.9, Section 5.2}{kida}]\label{thm-alternative-irreducible-amenable-non-amenable}
In Theorem \ref{thm-alternative-trivial-irreducible}, there exist two essentially unique invariant Borel maps $\varphi_{ia}, \varphi_{in}\colon Y\rightarrow \mathcal{F}_{0}(D)$ for $\mathcal{S}$ satisfying the following:
\begin{enumerate}
\item[(i)] $\varphi_{i}(x)=\varphi_{ia}(x)\cup \varphi_{in}(x)$ and $\varphi_{ia}(x)\cap \varphi_{in}(x)=\emptyset$ for a.e.\ $x\in Y$;
\item[(ii)] if $Q$ is in $F\in \mathcal{F}_{0}(D)$ with $\mu(\varphi_{ia}^{-1}(F))>0$, then 
\begin{enumerate}
\item[(a)] if $A$ is a Borel subset of $\varphi_{ia}^{-1}(F)$ with positive measure and $\psi \colon A\rightarrow M(\mathcal{PMF}(Q))$ is a $\rho_{Q}$-invariant Borel map for $\mathcal{S}$, then we see that $\psi(x)(\mathcal{MIN}(Q))=1$ for a.e.\ $x\in A$; 
\item[(b)] we have an essentially unique $\rho_{Q}$-invariant Borel map $\psi_{0}\colon \varphi_{ia}^{-1}(F)\rightarrow \partial_{2}C(Q)$ for $\mathcal{S}$ such that if $A$ is a Borel subset of $\varphi_{ia}^{-1}(F)$ with positive measure and $\psi \colon A\rightarrow M(\partial C(Q))$ is a $\rho_{Q}$-invariant Borel map for $\mathcal{S}$, then 
\[{\rm supp}(\psi(x))\subset {\rm supp}(\psi_{0}(x))\]
for a.e.\ $x\in A$;
\end{enumerate}
\item[(iii)] if $Q$ is in $F\in \mathcal{F}_{0}(D)$ with $\mu(\varphi_{in}^{-1}(F))>0$, then $\mathcal{S}$ admits neither $\rho_{Q}$-invariant Borel maps $A\rightarrow M(\mathcal{PMF}(Q))$ nor $A\rightarrow \partial_{2}C(Q)$ for any Borel subset $A\subset \varphi_{in}^{-1}(F)$ with positive measure.
\end{enumerate}
\end{thm}

We call $\varphi_{ia}$ (resp. $\varphi_{in}$) the system of irreducible and amenable or IA (resp. irreducible and non-amenable or IN) subsurfaces for $\mathcal{S}$. We often call an element in $\varphi_{ia}(x)$ (resp. $\varphi_{in}(x)$) an irreducible and amenable (resp. irreducible and non-amenable) subsurface for $x\in Y$ and in short, an IA (resp. IN) subsurface (or component). It is easy to see that if $A$ is a Borel subset of $Y$ with positive measure, then the IA, IN systems for $(\mathcal{S})_{A}$ are the restrictions of $\varphi_{ia}$, $\varphi_{in}$ to $A$, respectively (see \cite[Lemma 5.10]{kida}).

We recall some properties of IA components in the following lemma:

\begin{lem}[\ci{Lemma 5.13}{kida}]\label{lem-IA-normal}
Under the assumption $(\star)$, let $Y\subset X$ be a Borel subset with positive measure and let $\mathcal{S}$ be a subgroupoid of $(\mathcal{G})_{Y}$ of infinite type. Let $\mathcal{H}$ be a subgroupoid of $(\mathcal{G})_{Y}$ with $\mathcal{S}\vartriangleleft \mathcal{H}$. Suppose that $\mathcal{S}$ is reducible (thus, so is $\mathcal{H}$) and all the CRS, T, IA, IN systems for $\mathcal{S}$ are constant. Let $Q$ be an IA component for $\mathcal{S}$ and $\psi_{0}\colon Y\rightarrow \partial_{2}C(Q)$ be the $\rho_{Q}$-invariant Borel map for $\mathcal{S}$ as in Theorem \ref{thm-alternative-irreducible-amenable-non-amenable} (ii) (b). Then
\begin{enumerate}
\item[(i)] $\psi_{0}$ is $\rho_{Q}$-invariant for $\mathcal{H}$.
\item[(ii)] if we denote by $\psi \colon Y\rightarrow S(M)$ the CRS for $\mathcal{H}$, then $\sigma \subset \psi(x)$ for a.e.\ $x\in Y$, where $\sigma \in S(M)$ is the CRS for $\mathcal{S}$.
\item[(iii)] if we denote by $\psi_{ia}\colon Y\rightarrow \mathcal{F}_{0}(D)$ the IA system for $\mathcal{H}$, then $Q\in \psi_{ia}(x)$ for a.e.\ $x\in Y$.
\end{enumerate} 
\end{lem}

The following proposition implies that if a reducible subgroupoid has no IN components, then it is amenable as a groupoid:

\begin{prop}[\ci{Proposition 5.18}{kida}]\label{prop-aME act}
Under the assumption $(\star)$, let $Y\subset X$ be a Borel subset with positive measure and let $\mathcal{S}$ be a subgroupoid of $(\mathcal{G})_{Y}$ of infinite type. Suppose that $\mathcal{S}$ is reducible and there exists $\sigma \in S(M)$ such that $\rho(\gamma)\sigma =\sigma$ for a.e.\ $\gamma \in \mathcal{S}$. Let $\{ Q_{i}\}_{i}$ be the set of all components of $M_{\sigma}$ which are not pairs of pants and $\rho_{\sigma}\colon \mathcal{S}\rightarrow \prod_{i}\Gamma(Q_{i})$ be the induced cocycle $\prod_{i}\rho_{Q_{i}}$. Moreover, we assume that there exists a $\rho_{\sigma}$-invariant Borel map 
\[\psi \colon Y\rightarrow \prod_{i}\partial_{2}C(Q_{i}).\]
Then the groupoid $\mathcal{S}$ is amenable.
\end{prop}

Suppose that $\Gamma$ is a subgroup of $\Gamma(M; m)$ with an integer $m\geq 3$ and that $\sigma \in S(M)$ is fixed by each element of $\Gamma$. Let $p_{\sigma}\colon \Gamma \rightarrow \prod_{Q}\Gamma(Q)$ be the product $\prod_{Q}p_{Q}$, where $Q$ runs through all components in $M_{\sigma}$. Note that the kernel of $p_{\sigma}$ is contained in the amenable subgroup of $\Gamma(M)$ generated by Dehn twists about all curves in $\sigma$ by Lemma \ref{lem-ker-dehn}. Thus, it is easily shown that if every component of $M_{\sigma}$ is either T or IA, then $\Gamma$ is amenable, which implies Proposition \ref{prop-aME act} in the case where $X$ is a point.


\subsection{Actions of hyperbolic groups}

In this subsection, we consider subgroupoids of a groupoid defined by a non-singular action of a hyperbolic group. In subsequent sections, we only use Lemma \ref{lem-free-prod}, in which we consider a groupoid arising from an action of a free group of rank $2$.

\begin{assumption}
We call the following assumption $(\star)_{\rm h}$: let $\Gamma$ be an infinite subgroup of a hyperbolic group $\Gamma_{0}$. Let $(X, \mu)$ be a standard finite measure space and assume that we have a non-singular $\Gamma$-action on $(X, \mu)$. We denote by $\mathcal{G}=\Gamma \ltimes X$ and $\rho \colon \mathcal{G}\rightarrow \Gamma$ the associated groupoid and cocycle, respectively. 
\end{assumption}

For a hyperbolic group $\Gamma_{0}$, let $\partial \Gamma_{0}$ be the boundary at infinity and let $M(\partial \Gamma_{0})$ be the space of all probability measures on $\partial \Gamma_{0}$. We denote by $\partial_{2}\Gamma_{0}$ the quotient space of $\partial \Gamma_{0}\times \partial \Gamma_{0}$ by the coordinate exchanging action of the symmetric group of two letters, which can naturally be viewed as a Borel subset of $M(\partial \Gamma_{0})$ as in the case of the boundary of the curve complex. We can show the following proposition with the methods from the proof of Theorem \ref{thm-ia}:

\begin{prop}\label{prop-hyp-action}
Under the assumption $(\star)_{\rm h}$, let $Y$ be a Borel subset of $X$ with positive measure and let $\mathcal{S}$ be a subgroupoid of $(\mathcal{G})_{Y}$ of infinite type. Assume that there is a $\rho$-invariant Borel map $Y\rightarrow M(\partial \Gamma_{0})$ for $\mathcal{S}$. Then
\begin{enumerate}
\item[(i)] there exists an essentially unique $\rho$-invariant Borel map $\varphi_{0}\colon Y\rightarrow \partial_{2}\Gamma_{0}$ for $\mathcal{S}$ satisfying the following: if $Y'$ is a Borel subset of $Y$ with positive measure and $\varphi \colon Y'\rightarrow M(\partial \Gamma_{0})$ is a $\rho$-invariant Borel map for $\mathcal{S}$, then
\[{\rm supp}(\varphi(x))\subset {\rm supp}(\varphi_{0}(x))\]
for a.e.\ $x\in Y'$.
\item[(ii)] if $\mathcal{T}$ is a subgroupoid of $(\mathcal{G})_{Y}$ with $\mathcal{S}\vartriangleleft \mathcal{T}$, then $\varphi_{0}$ is invariant also for $\mathcal{T}$.
\item[(iii)] $\mathcal{S}$ is amenable as a groupoid.  
\end{enumerate}
\end{prop}

Using Lemma \ref{lem-iso-normal} and Proposition \ref{prop-hyp-action} (ii), (iii), we can show the following as in Corollary \ref{cor-mcg-stab}. Note that the set consisting of all points in $\partial_{2}\Gamma_{0}$ fixed by some infinite subgroup of $\Gamma_{0}$ is countable (see \cite[Chapitre 8]{ghys-harpe} or the comment about the dynamics of a hyperbolic group on its boundary right after \cite[Theorem 3.3]{kida}).

\begin{cor}\label{cor-hyp-stab}
Under the assumption $(\star)_{\rm h}$, suppose that $\Gamma$ is non-amenable and that the $\Gamma$-action on $(X, \mu)$ is measure-preserving. Then the isotropy group 
\[\mathcal{G}_{x}^{x}=\{ \gamma \in \mathcal{G}: r(\gamma)=s(\gamma)=x\}\]
is either finite or non-amenable for a.e.\ $x\in X$.
\end{cor}

\begin{lem}\label{lem-free-prod}
Let $G_{1}$, $G_{2}$ be two infinite cyclic groups and suppose that we have a measure-preserving action of the free product $G=G_{1}*G_{2}$ on a standard finite measure space $(X, \mu)$. Let $\mathcal{G}$, $\mathcal{G}_{1}$ and $\mathcal{G}_{2}$ be the groupoids arising from the actions of $G$, $G_{1}$ and $G_{2}$, respectively. Then the subgroupoid $(\mathcal{G}_{1})_{A}\vee (\mathcal{G}_{2})_{A}$ of $(\mathcal{G})_{A}$ generated by $(\mathcal{G}_{1})_{A}$ and $(\mathcal{G}_{2})_{A}$ is non-amenable for any Borel subset $A\subset X$ with positive measure.
\end{lem}

\begin{proof}
Suppose that $(\mathcal{G}_{1})_{A}\vee (\mathcal{G}_{2})_{A}$ is amenable. We have the natural cocycle $\rho \colon \mathcal{G}\rightarrow G$. It follows that we can find a $\rho$-invariant Borel map $\varphi_{0}\colon A\rightarrow \partial_{2}G$ for $(\mathcal{G}_{1})_{A}\vee (\mathcal{G}_{2})_{A}$ as in Proposition \ref{prop-hyp-action} (i). Let $a_{i}^{\pm}\in \partial G$ be the two fixed points on the boundary $\partial G$ of $G$ for the action of the group $G_{i}$ for $i=1, 2$. Then the constant map $\varphi_{i}\colon A\rightarrow \partial_{2}G$ with value $\{ a_{i}^{\pm}\}$ is $\rho$-invariant for the subgroupoid $(\mathcal{G}_{i})_{A}$ of infinite type. It follows that $\varphi_{i}$ has to satisfy the property in Proposition \ref{prop-hyp-action} (i). Thus, we have ${\rm supp}(\varphi(x))\subset {\rm supp}(\varphi_{i}(x))=\{ a_{i}^{\pm}\}$ for $i=1, 2$, which is a contradiction because $\{ a_{1}^{\pm}\} \cap \{ a_{2}^{\pm}\} =\emptyset$. 
\end{proof}


\section{Characterizations of reducible subgroupoids}\label{sec-cha-red}

The next two propositions characterize an amenable and a non-amenable reducible subgroupoid, respectively, in terms of amenability, non-amenability and normality. As in the previous section, we use the following notation under the assumption $(\star)$:
\begin{itemize}
\item For a subgroup $\Gamma'$ of $\Gamma$, let $\mathcal{G}_{\Gamma'}$ denote the subgroupoid of $\mathcal{G}$ generated by the action of $\Gamma'$.  
\item For $\sigma \in S(M)$, we denote by $D_{\sigma}$ the intersection of $\Gamma$ and the subgroup generated by Dehn twists about all curves in $\sigma$. We write $\mathcal{G}_{\sigma}$ instead of $\mathcal{G}_{D_{\sigma}}$ for simplicity. If $\sigma$ consists of one element $\alpha \in V(C)$, then we write $D_{\alpha}$ (resp. $\mathcal{G}_{\alpha}$) instead of $D_{\sigma}$ (resp. $\mathcal{G}_{\sigma}$).
\item For $\sigma \in S(M)$, we write
\[\Gamma_{\sigma}=\{ g\in \Gamma : g\sigma =\sigma\}.\]
\end{itemize}

\begin{prop}\label{prop-cha-red-ame}
Under the assumption $(\star)$, let $Y\subset X$ be a Borel subset with positive measure and let $\mathcal{S}$ be a subgroupoid of $(\mathcal{G})_{Y}$ of infinite type. Suppose that $\mathcal{S}$ is amenable. Consider the following two assertions:
\begin{enumerate}
\item[(i)] $\mathcal{S}$ is reducible.
\item[(ii)] For any Borel subset $A$ of $Y$ with positive measure, we have a Borel subset $B$ of $A$ with positive measure and the following three subgroupoids $\mathcal{S}'$, $\mathcal{S}''$ and $\mathcal{T}$ of $(\mathcal{G})_{B}$: 
\begin{enumerate}
\item[(a)] an amenable subgroupoid $\mathcal{S}'$ with $(\mathcal{S})_{B}<\mathcal{S}'$;
\item[(b)] a subgroupoid $\mathcal{S}''$ of infinite type with $\mathcal{S}''<\mathcal{S}'$;
\item[(c)] a non-amenable subgroupoid $\mathcal{T}$ with $\mathcal{S}''\vartriangleleft \mathcal{T}$. 
\end{enumerate}
\end{enumerate}
Then the assertion {\rm (ii)} implies the assertion {\rm (i)}. If $\Gamma$ is a subgroup of finite index in $\Gamma(M;m)$ and the $\Gamma$-action on $(X, \mu)$ is measure-preserving, then the converse also holds.
\end{prop}

\begin{proof}
First, we show that the assertion (ii) implies the assertion (i). If $\mathcal{S}$ were not reducible, then since $\mathcal{S}$ is of infinite type and amenable, there would exist an invariant Borel subset $A$ of $Y$ for $\mathcal{S}$ with positive measure such that $(\mathcal{S})_{A}$ is IA (see Theorem \ref{thm-alternative}). It follows from the assumption (ii) that we have a Borel subset $B$ of $A$ with positive measure and subgroupoids $\mathcal{S}'$, $\mathcal{S}''$ and $\mathcal{T}$ satisfying the conditions in (ii). Since $(\mathcal{S})_{A}$ is IA and $\mathcal{S}'$ is amenable, it follows from Theorem \ref{thm-alternative} that $\mathcal{S}'$ is IA. Moreover, $\mathcal{S}''$ is also IA. Thus, $\mathcal{T}$ is also IA and amenable by Theorem \ref{thm-ia} (ii) and Proposition \ref{prop-bou-IA}, which is a contradiction. 

Next, we assume that $\Gamma$ is a subgroup of finite index in $\Gamma(M;m)$ and show that the assertion (i) implies the assertion (ii). Let $A$ be a Borel subset of $Y$ with positive measure. Then there exists a Borel subset $B$ of $A$ with positive measure satisfying the following conditions (see Lemma \ref{lem-T-strong} for the second condition):
\begin{itemize}
\item all of the CRS, T and IA systems for $\mathcal{S}$ are constant on $B$. Let $\sigma \in S(M)$, $\varphi_{t}, \varphi_{ia}\in {\mathcal F}_{0}(M)$ be their values on $B$, respectively. Note that the IN system for $\mathcal{S}$ is empty since $\mathcal{S}$ is amenable;
\item for a.e.\ $\gamma \in (\mathcal{S})_{B}$ and any component $Q$ in $\varphi_{t}$ and $\alpha \in V(C(Q))$, we have $\rho_{Q}(\gamma)\alpha =\alpha$, where $\rho_{Q}\colon (\mathcal{S})_{B}\rightarrow \Gamma(Q)$ is the composition of $\rho$ and the natural projection $\Gamma_{\sigma}\rightarrow \Gamma(Q)$. 
\end{itemize}

For each $Q\in \varphi_{ia}$, we have the canonical invariant Borel map $\psi_{Q}\colon B\rightarrow \partial_{2}C(Q)$ for $(\mathcal{S})_{B}$ as in Theorem \ref{thm-alternative-irreducible-amenable-non-amenable} (ii) (b). Let $\mathcal{S}'$ be a subgroupoid of $(\mathcal{G})_{B}$ consisting of all $\gamma \in (\mathcal{G})_{B}$ satisfying 
\[\rho(\gamma)\sigma =\sigma,\ \ \rho_{Q}(\gamma)\psi_{Q}(s(\gamma))=\psi_{Q}(r(\gamma)),\ \ \rho(\gamma)\alpha =\alpha\]
for any $Q\in \varphi_{ia}$, any $\alpha \in V(C(R))$ and any $R\in \varphi_{t}$ which is not a pair of pants. Note that $(\mathcal{S})_{B}<\mathcal{S}'$. It follows from Proposition \ref{prop-aME act} that $\mathcal{S}'$ is amenable.

If $|\sigma|<\kappa(M)+1$, then let $\mathcal{S}''=(\mathcal{G}_{\sigma})_{B}$.  Then $\mathcal{S}''<\mathcal{S}'$. Since $\Gamma$ is a subgroup of finite index in $\Gamma(M;m)$ and there exists a component of $M_{\sigma}$ which is not a pair of pants, we see that $\mathcal{S}''$ is of infinite type and $\Gamma_{\sigma}$ is non-amenable. Thus, the subgroupoid $\mathcal{T}=(\mathcal{G}_{\Gamma_{\sigma}})_{B}$ is non-amenable. Moreover, $\mathcal{S}''\vartriangleleft \mathcal{T}$ since $D_{\sigma}$ is a normal subgroup of $\Gamma_{\sigma}$ by Lemma \ref{lem-ker-dehn}. This completes the construction of subgroupoids in the assertion (ii) in the case of $|\sigma|<\kappa(M)+1$.

If $|\sigma|=\kappa(M)+1$, then $\mathcal{S}'=(\mathcal{G}_{\sigma})_{B}$ and it is amenable. Choose $\alpha_{0}\in \sigma$. Let $\sigma'=\sigma \setminus \{ \alpha_{0}\}$, which is an element in $S(M)$ since $\kappa(M)>0$. Then $\mathcal{S}''=(\mathcal{G}_{\sigma'})_{B}$ is a subgroupoid of infinite type with $\mathcal{S}''<\mathcal{S}'$. Define $\mathcal{T}=(\mathcal{G}_{\Gamma_{\sigma'}})_{B}$. Then $\mathcal{T}$ is non-amenable and $\mathcal{S}''\vartriangleleft \mathcal{T}$ since $D_{\sigma'}$ is a normal subgroup of $\Gamma_{\sigma'}$ by Lemma \ref{lem-ker-dehn}. This completes the construction of subgroupoids in the assertion (ii) in the case of $|\sigma|=\kappa(M)+1$. 
\end{proof}

\begin{prop}\label{prop-cha-red-non}
Under the assumption $(\star)$, let $Y\subset X$ be a Borel subset with positive measure and let $\mathcal{S}$ be a subgroupoid of $(\mathcal{G})_{Y}$ of infinite type. Suppose that $(\mathcal{S})_{Y'}$ is not amenable for any Borel subset $Y'$ of $Y$ with positive measure. Consider the following two assertions:
\begin{enumerate}
\item[(i)] $\mathcal{S}$ is reducible.
\item[(ii)] For any Borel subset $A$ of $Y$ with positive measure, we have a Borel subset $B$ of $A$ with positive measure and the following two subgroupoids $\mathcal{S}'$ and $\mathcal{S}''$ of $(\mathcal{G})_{B}$: 
\begin{enumerate}
\item[(a)] a subgroupoid $\mathcal{S}'$ with $(\mathcal{S})_{B}<\mathcal{S}'$;
\item[(b)] an amenable subgroupoid $\mathcal{S}''$ of infinite type with $\mathcal{S}''\vartriangleleft \mathcal{S}'$.
\end{enumerate}
\end{enumerate}
Then the assertion {\rm (ii)} implies the assertion {\rm (i)}. If $\Gamma$ is a subgroup of finite index in $\Gamma(M;m)$ and the $\Gamma$-action on $(X, \mu)$ is measure-preserving, then the converse also holds.
\end{prop}

\begin{proof}
First, we show that the assertion (ii) implies the assertion (i). Suppose that $\mathcal{S}$ is not reducible. Then there exists a Borel subset $A$ of $Y$ with positive measure such that for any Borel subset $B$ of $A$ with positive measure, there is no invariant Borel map $B\rightarrow M(\mathcal{PMF})$ for $\mathcal{S}$ (see Remark \ref{rem-div}). By the assumption (ii), we have a Borel subset $B$ of $A$ with positive measure and two subgroupoids $\mathcal{S}'$ and $\mathcal{S}''$ satisfying the conditions in (ii). Since $\mathcal{S}''$ is amenable, by Theorem \ref{thm-alternative}, we have a Borel partition $B=B_{1}\sqcup B_{2}$ (up to null sets) such that $(\mathcal{S}'')_{B_{1}}$ is IA and $(\mathcal{S}'')_{B_{2}}$ is reducible. It follows from Theorem \ref{thm-ia} (ii) and Theorem \ref{thm-reducible} (iii) that $(\mathcal{S}')_{B_{1}}$ is IA and $(\mathcal{S}')_{B_{2}}$ is reducible. If $B_{1}$ has positive measure, then $(\mathcal{S})_{B_{1}}$ is non-amenable by the assumption of $\mathcal{S}$. Since $(\mathcal{S})_{B}<\mathcal{S}'$, the relation $(\mathcal{S}')_{B_{1}}$ is non-amenable, which contradicts to Theorem \ref{thm-ia} (iii). On the other hand, if $B_{2}$ has positive measure, then $(\mathcal{S})_{B_{2}}$ has an invariant Borel map $B_{2}\rightarrow S(M)\subset M(\mathcal{PMF})$, which is also a contradiction.

Next, we assume that $\Gamma$ is a subgroup of finite index in $\Gamma(M;m)$ and show that the converse also holds. Let $A$ be a Borel subset of $Y$ with positive measure. Then there exists a Borel subset $B$ of $A$ with positive measure such that the CRS for $\mathcal{S}$ is constant on $B$. We denote by $\sigma \in S(M)$ its value on $B$. Define a subgroupoid
\[\mathcal{S}'=\{\gamma \in (\mathcal{G})_{B}: \rho(\gamma)\sigma =\sigma \} =(\mathcal{G}_{\Gamma_{\sigma}})_{B},\]
which satisfies $(\mathcal{S})_{B}<\mathcal{S}'$. Let $\mathcal{S}''=(\mathcal{G}_{\sigma})_{B}$. Then $\mathcal{S}''$ is of infinite type since $\Gamma$ is a subgroup of finite index in $\Gamma(M; m)$ and the $\Gamma$-action on $(X, \mu)$ is measure-preserving. Since $D_{\sigma}$ is a normal subgroup of $\Gamma_{\sigma}$ and it is amenable by Lemma \ref{lem-ker-dehn}, we see that $\mathcal{S}''\vartriangleleft \mathcal{S}'$ and $\mathcal{S}''$ is amenable. 
\end{proof}

\begin{assumption}\label{assumption-bullet0}
We call the following assumption $(\bullet)$: for $i=1, 2$, let $\Gamma_{i}$ be a finite index subgroup of $\Gamma(M_{i}; m_{i})$, where $M_{i}$ is a surface with $\kappa(M_{i})>0$ and $m_{i}\geq 3$ is an integer. Consider a measure-preserving action of $\Gamma_{i}$ on a standard finite measure space $(X_{i}, \mu_{i})$ and let 
\[\mathcal{G}^{i}=\mathcal{G}^{i}_{\Gamma},\ \ \rho_{i} \colon \mathcal{G}^{i}\rightarrow \Gamma_{i}\]
be the induced groupoid, cocycle, respectively. Suppose that we have a groupoid isomorphism 
\[f\colon (\mathcal{G}^{1})_{Y_{1}}\rightarrow (\mathcal{G}^{2})_{Y_{2}},\]
where $Y_{i}\subset X_{i}$ is a Borel subset satisfying $\mathcal{G}^{i}Y_{i}=X_{i}$ up to null sets.
\end{assumption}

The following corollary is a consequence of Propositions \ref{prop-cha-red-ame} and \ref{prop-cha-red-non} which characterize reducible subrelations:

\begin{cor}\label{cor-red-pre}
Under the assumption $(\bullet)$, let $A_{1}$ be a Borel subset of $Y_{1}$ with positive measure and let $\mathcal{S}^{1}$ be a subgroupoid of $(\mathcal{G}^{1})_{A_{1}}$ of infinite type. Then $\mathcal{S}^{1}$ is reducible if and only if the image $f(\mathcal{S}^{1})$ is reducible.
\end{cor}

Next, we characterize maximal reducible subgroupoids. In the assumption $(\star)$, let $Y$ be a Borel subset of $X$ with positive measure and let $\varphi \colon Y\rightarrow S(M)$ be a Borel map. Then we define a reducible subgroupoid
\[\mathcal{S}_{\varphi}=\{ \gamma \in (\mathcal{G})_{Y}: \rho(\gamma)\varphi(s(\gamma))=\varphi(r(\gamma))\}.\]

\begin{prop}\label{prop-gen-red-crs}
Under the assumption $(\star)$, let $Y$ be a Borel subset of $X$ with positive measure and let $\varphi \colon Y\rightarrow S(M)$ be a Borel map. Assume that $\Gamma$ is a subgroup of finite index in $\Gamma(M; m)$ and that the $\Gamma$-action on $(X, \mu)$ is measure-preserving. Then the CRS for $\mathcal{S}_{\varphi}$ is $\varphi$ and for a.e.\ $x\in Y$, each component of $M_{\varphi(x)}$ either is a pair of pants or is IN for $\mathcal{S}_{\varphi}$.  
\end{prop}

\begin{proof}
We may assume that all of the CRS, T, IA and IN systems for $\mathcal{S}_{\varphi}$ and $\varphi$ are constant by \cite[Lemmas 5.7, 5.10]{kida}. We denote the value of $\varphi$ by the same symbol. Note that $\mathcal{S}_{\varphi}$ is equal to $(\mathcal{G}_{\Gamma_{\varphi}})_{Y}$. It follows from Lemmas \ref{lem-red-crs-group}, \ref{lem-red-crs} that the CRS for $\mathcal{S}_{\varphi}$ is $\varphi$. 

Let $Q$ be a component of $M_{\varphi}$ which is not a pair of pants. Let $g_{1}, g_{2}\in \Gamma_{\varphi}$ be elements such that $\{ p_{Q}(g_{1}), p_{Q}(g_{2})\}$ is an independent pair of pseudo-Anosov elements in $\Gamma(Q)$, where $p_{Q}\colon \Gamma_{\varphi}\rightarrow \Gamma(Q)$ is the natural homomorphism. Let $G$ be the subgroup of $\Gamma_{\varphi}$ generated by $g_{1}$ and $g_{2}$. Note that $(\mathcal{G}_{G})_{Y}<\mathcal{S}_{\varphi}$. 

If $Q$ were T for $\mathcal{S}_{\varphi}$, then it follows from Lemma \ref{lem-T-strong} that there would exist a Borel subset $A$ of $Y$ with positive measure such that $\rho_{Q}(\gamma)\alpha =\alpha$ for any $\alpha \in V(C(Q))$ and for a.e.\ $\gamma \in (\mathcal{S}_{\varphi})_{A}$, where $\rho_{Q}$ is the composition of $\rho$ and $p_{Q}$. This contradicts the fact that $p_{Q}(g_{1}^{n})\alpha \neq \alpha$ for any $\alpha \in V(C(Q))$ and all $n\in \mathbb{Z}\setminus \{ 0\}$. 

If $Q$ were IA for $\mathcal{S}_{\varphi}$, then we would have the canonical invariant Borel map $\phi \colon Y\rightarrow \partial_{2}C(Q)$ for $\mathcal{S}_{\varphi}$ as in Theorem \ref{thm-alternative-irreducible-amenable-non-amenable} (ii) (b). For $i=1, 2$, define a Borel map $\phi_{i}\colon Y\rightarrow \partial_{2}C(Q)$ to be the constant map whose value is the image of $\{ F_{\pm}(p_{Q}(g_{i}))\}$ in $\partial C(Q)$. Recall that the natural map $\mathcal{MIN}\rightarrow \partial C$ is injective on the set of all pseudo-Anosov foliations. It follows that $\phi_{i}$ is the canonical invariant Borel map for $(\mathcal{G}_{G_{i}})_{Y}$ for $i=1, 2$, where $G_{i}$ is the cyclic subgroup generated by $g_{i}$. Since $\phi$ is invariant for $(\mathcal{G}_{G_{i}})_{Y}$, we have the inclusion
\[{\rm supp}(\phi(x))\subset {\rm supp}(\phi_{i}(x))\]
for a.e.\ $x\in Y$ and any $i=1, 2$. This is a contradiction because $\{ F_{\pm}(p_{Q}(g_{1}))\} \cap \{ F_{\pm}(p_{Q}(g_{2}))\} =\emptyset$.   
\end{proof}

In what follows, we regard $V(C)$ as a subset of $S(M)$ naturally.

\begin{cor}\label{cor-max-red}
Under the assumption $(\star)$, let $Y$ be a Borel subset of $X$ with positive measure and let $\varphi \colon Y\rightarrow V(C)$ be a Borel map. Assume that $\Gamma$ is a finite index subgroup of $\Gamma(M; m)$ and that the $\Gamma$-action on $(X, \mu)$ is measure-preserving. If $\mathcal{S}$ is a reducible subgroupoid of $(\mathcal{G})_{Y}$ with $\mathcal{S}_{\varphi}<\mathcal{S}$, then $\mathcal{S}=\mathcal{S}_{\varphi}$.
\end{cor}

\begin{proof}
Let $\psi \colon Y\rightarrow S(M)$ be the CRS for $\mathcal{S}$. It is enough to show $\varphi =\psi$ up to null sets. Choose $\alpha \in V(C)$ and $\sigma \in S(M)$ such that $\mu(\varphi^{-1}(\alpha)\cap \psi^{-1}(\sigma))>0$ and put $A=\varphi^{-1}(\alpha)\cap \psi^{-1}(\sigma)$. It suffices to prove $\varphi =\psi$ a.e.\ on $A$, that is, $\sigma =\{ \alpha \}$. We may assume that all the T, IA and IN systems for $\mathcal{S}_{\varphi}$ on $A$ are constant.

Choose $\beta \in \sigma$. Since $\beta$ is in the CRS for $(\mathcal{S})_{A}$, the pair $(\beta, A)$ is $\rho$-invariant for $\mathcal{S}_{\varphi}$. If we had a component $Q$ of $M_{\alpha}$ which is not a pair of pants and satisfies $r(\beta, Q)\neq \emptyset$, then $Q$ would be T for $\mathcal{S}_{\varphi}$ by Theorem \ref{thm-alternative-trivial-irreducible}, which contradicts Proposition \ref{prop-gen-red-crs}. Thus, $r(\beta, Q)=\emptyset$ for each component $Q$ of $M_{\alpha}$ which is not a pair of pants. It follows from Lemma \ref{lem-pants} that $\beta$ is a boundary component of $Q$, that is, $\alpha =\beta$. Therefore, $\sigma =\{ \alpha \}$ and $\varphi =\psi$ a.e.\ on $A$.  
\end{proof}

\begin{lem}\label{lem-max-map}
Under the assumption $(\star)$, let $Y\subset X$ be a Borel subset with positive measure and let $\mathcal{S}$ be a subgroupoid of $(\mathcal{G})_{Y}$ of infinite type. Suppose that $\mathcal{S}$ is reducible. Then there exists a Borel map $\psi \colon Y\rightarrow V(C)$ such that $\mathcal{S}<\mathcal{S}_{\psi}$.
\end{lem}

\begin{proof}
Let $\varphi \colon Y\rightarrow S(M)$ be the  CRS for $\mathcal{S}$. Choose a countable Borel partition $Y=\bigsqcup Y_{n}$ of $Y$ such that $\varphi$ is constant on each $Y_{n}$. Let $\alpha_{n}\in V(C)$ be an element such that $\alpha_{n}\in \varphi(x)$ for a.e.\ $x\in Y_{n}$. Then the constant map $Y_{n}\ni x\mapsto \alpha_{n}\in V(C)$ is invariant for $(\mathcal{S})_{Y_{n}}$. It follows from the proof of Lemma \ref{lem-ext-inv} that we can construct an invariant Borel map $\psi \colon Y\rightarrow V(C)$ for $\mathcal{S}$.
\end{proof}

The following corollary is a consequence of Corollaries \ref{cor-red-pre}, \ref{cor-max-red} and Lemma \ref{lem-max-map}:

\begin{cor}\label{cor-pre-max-map}
Under the assumption $(\bullet)$, let $A_{1}$ be a Borel subset of $Y_{1}$ with positive measure and let $\varphi_{1} \colon A_{1}\rightarrow V(C(M_{1}))$ be a Borel map. Put $A_{2}=f(A_{1})$ and 
\[\mathcal{S}_{\varphi_{1}}^{1}=\{ \gamma \in (\mathcal{G}^{1})_{Y_{1}}: \rho_{1}(\gamma)\varphi_{1}(s(\gamma))=\varphi_{1}(r(\gamma))\}.\]
Then there exists a Borel map $\varphi_{2} \colon A_{2}\rightarrow V(C(M_{2}))$ such that $f(\mathcal{S}_{\varphi_{1}}^{1})=\mathcal{S}_{\varphi_{2}}^{2}$, where
\[\mathcal{S}_{\varphi_{2}}^{2}=\{ \gamma \in (\mathcal{G}^{2})_{Y_{2}}: \rho_{2}(\gamma)\varphi_{2}(s(\gamma))=\varphi_{2}(r(\gamma))\}.\]
\end{cor}


\section{An equivariant Borel map from a self ME coupling}\label{sec-eq-ME}

In the next lemma, we study a normal amenable subgroupoid of a maximal reducible subgroupoid. As in the previous section, we regard the vertex set $V(C)$ as a subset of the simplex set $S(M)$ naturally.

\begin{lem}\label{lem-dehn-max-normal}
Under the assumption $(\star)$, let $Y$ be a Borel subset of $X$ with positive measure and let $\varphi \colon Y\rightarrow V(C)$ be a Borel map. Assume that $\Gamma$ is a finite index subgroup of $\Gamma(M; m)$ and that the $\Gamma$-action on $(X, \mu)$ is measure-preserving. If $\mathcal{S}$ is an amenable subgroupoid of $\mathcal{S}_{\varphi}$ of infinite type with $\mathcal{S}\vartriangleleft \mathcal{S}_{\varphi}$, then there exists a countable Borel partition $Y=\bigsqcup Y_{n}$ of $Y$ satisfying the following conditions:
\begin{enumerate}
\item[(i)] the map $\varphi$ is constant a.e.\ on $Y_{n}$. Let $\alpha_{n}\in V(C)$ be its value;
\item[(ii)] for each $n$, we have $(\mathcal{S})_{Y_{n}}<(\mathcal{G}_{\alpha_{n}})_{Y_{n}}<(\mathcal{S}_{\varphi})_{Y_{n}}$.
\end{enumerate}
\end{lem} 
 
\begin{proof}
Recall that $\mathcal{S}_{\varphi}$ is reducible and its CRS is given by $\varphi$ (see Proposition \ref{prop-gen-red-crs}). Since $\mathcal{S}$ is a subgroupoid of $\mathcal{S}_{\varphi}$, it is also reducible. Let $\psi \colon Y\rightarrow S(M)$ be the CRS for $\mathcal{S}$. Since $\mathcal{S}$ is normal in $\mathcal{S}_{\varphi}$, the map $\psi$ is invariant also for $\mathcal{S}_{\varphi}$ by Theorem \ref{thm-reducible} (iii). By the definition of essential $\rho$-invariant pairs, we see that $\psi(x)\subset \varphi(x)$ for a.e.\ $x\in Y$. It follows that $\psi(x)=\varphi(x)$ for a.e.\ $x\in Y$. 

Let $A$ be a Borel subset of $Y$ with positive measure such that all of the CRS $\varphi =\psi$ and T, IA and IN systems for $\mathcal{S}$ and $\mathcal{S}_{\varphi}$ are constant on $A$. We denote by $\alpha \in V(C)$ the value of $\varphi =\psi$ on $A$. If $Q$ is a component of $M_{\alpha}$, then $Q$ is not IN for $(\mathcal{S})_{A}$ since $\mathcal{S}$ is amenable. If $Q$ were IA for $(\mathcal{S})_{A}$, then $Q$ would be IA also for $(\mathcal{S}_{\varphi})_{A}$ by Lemma \ref{lem-IA-normal} (iii), which contradicts Proposition \ref{prop-gen-red-crs}. Thus, each component of $M_{\alpha}$ is T for $(\mathcal{S})_{A}$.

It follows from Lemma \ref{lem-T-strong} that we have a countable Borel partition $A=\bigsqcup A_{n}$ of $A$ such that $\rho(\gamma)\beta =\beta$ for each component $Q$ of $M_{\alpha}$ and $\beta \in V(C(Q))$ and for a.e.\ $\gamma \in (\mathcal{S})_{A_{n}}$ for any $n$. For a.e.\ $\gamma \in (\mathcal{S})_{A_{n}}$, consider the subgroup of $\Gamma$ generated by $\rho(\gamma)$. If $\rho(\gamma)$ is non-trivial, then the CRS for the subgroup is $\{ \alpha \}$ by Lemma \ref{lem-1-crs}. It follows from \cite[Corollary 7.18]{ivanov1} (see the comment right before Lemma \ref{lem-red-crs-group}) that $\rho(\gamma)$ lies in the kernel of the natural homomorphism from $\Gamma_{\alpha}$ into $\prod_{Q}\Gamma(Q)$, where $Q$ is taken over all components of $M_{\alpha}$. Thus, $\rho(\gamma)\in D_{\alpha}$ by Lemma \ref{lem-ker-dehn}. Since $A$ is any Borel subset of $Y$ with positive measure such that all of $\varphi =\psi$ and T, IA and IN systems for $\mathcal{S}$ and $\mathcal{S}_{\varphi}$ are constant on $A$, we complete the proof.
\end{proof}

Under the assumption $(\bullet)$, let $\alpha \in V(C(M_{1}))$. Define the constant map $\varphi_{\alpha}\colon Y_{1}\ni x\mapsto \alpha \in V(C(M_{1}))$. It follows from Corollary \ref{cor-pre-max-map} that we have a Borel map $\varphi_{2} \colon Y_{2}\rightarrow V(C(M_{2}))$ such that $f(\mathcal{S}_{\varphi_{\alpha}}^{1})=\mathcal{S}_{\varphi_{2}}^{2}$, where we use the same notation as in the corollary. Since the intersection of $\Gamma_{1}$ and the subgroup of $\Gamma(M_{1})$ generated by a Dehn twist along $\alpha$ is normal in 
\[\Gamma_{1, \alpha}=\{ g\in \Gamma_{1}: g\alpha =\alpha \} \] 
by Lemma \ref{lem-ker-dehn} and $\mathcal{S}_{\varphi_{\alpha}}^{1}=(\mathcal{G}_{\Gamma_{1, \alpha}}^{1})_{Y_{1}}$, we see that $(\mathcal{G}_{\alpha}^{1})_{Y_{1}}\vartriangleleft \mathcal{S}_{\varphi_{\alpha}}^{1}$. Thus, $f((\mathcal{G}_{\alpha}^{1})_{Y_{1}})\vartriangleleft \mathcal{S}_{\varphi_{2}}^{2}$. By Lemma \ref{lem-dehn-max-normal}, we have a countable Borel partition $Y_{2}=\bigsqcup A_{n}$ such that
\begin{enumerate}
\item[(i)] the map $\varphi_{2}$ is constant on $A_{n}$ for each $n$. Let us denote the value by $\beta_{n}\in V(C(M_{2}))$;
\item[(ii)] for each $n$, we have $(f((\mathcal{G}_{\alpha}^{1})_{Y_{1}}))_{A_{n}}<(\mathcal{G}_{\beta_{n}}^{2})_{A_{n}}<(\mathcal{S}_{\varphi_{2}}^{2})_{A_{n}}$.
\end{enumerate}
Therefore, for each $\alpha \in V(C(M_{1}))$, we can define a Borel map
\[\Psi(\cdot, \alpha)\colon Y_{1}\rightarrow V(C(M_{2}))\]
by $\Psi(x, \alpha)=\beta_{n}$ if $x\in f^{-1}(A_{n})$ (up to null sets). Note that this map does not depend on the decomposition $Y_{2}=\bigsqcup A_{n}$.

\begin{lem}\label{lem-Psi1}
If $\alpha, \alpha'\in V(C(M_{1}))$ satisfy $i(\alpha, \alpha')=0$, then $i(\Psi(x, \alpha), \Psi(x, \alpha'))=0$ for a.e.\ $x\in Y_{1}$.
\end{lem}
 
\begin{proof}
Since $i(\alpha, \alpha')=0$, we see that 
\[(\mathcal{G}_{\alpha}^{1})_{A} \vartriangleleft (\mathcal{G}_{\alpha}^{1})_{A}\vee (\mathcal{G}_{\alpha'}^{1})_{A}\]
for any Borel subset $A$ of $Y_{1}$ with positive measure (see Lemma \ref{lem-normal-dir-prod}). It follows from the construction of $\Psi(\cdot, {\alpha})$, $\Psi(\cdot, {\alpha'})$ that we have a countable Borel partition $Y_{2}=\bigsqcup A_{n}$ and $\beta_{n}, \beta_{n}'\in V(C(M_{2}))$ such that
\[(f((\mathcal{G}_{\alpha}^{1})_{Y_{1}}))_{A_{n}}<(\mathcal{G}_{\beta_{n}}^{2})_{A_{n}},\ \ \ (f((\mathcal{G}_{\alpha'}^{1})_{Y_{1}}))_{A_{n}}<(\mathcal{G}_{\beta_{n}'}^{2})_{A_{n}}\]
for each $n$. 
Using Lemma \ref{lem-crs-dehn}, we see that $(f((\mathcal{G}_{\alpha}^{1})_{Y_{1}}))_{A_{n}}$ (resp. $(f((\mathcal{G}_{\alpha'}^{1})_{Y_{1}}))_{A_{n}}$) is a reducible subrelation of $(\mathcal{G}^{2})_{A_{n}}$ and its CRS is given by the constant map $A_{n}\ni x\mapsto \beta_{n}$ (resp. $\beta_{n}'$) $\in V(C(M_{2}))$. It follows from the above normality that the constant map $A_{n}\ni x\mapsto \beta_{n}$ is invariant also for $(f((\mathcal{G}_{\alpha'}^{1})_{Y_{1}}))_{A_{n}}$, which implies $i(\beta_{n}, \beta_{n}')=0$ by the pureness of the pair $(\beta_{n}', A_{n})$ for $(f((\mathcal{G}_{\alpha'}^{1})_{Y_{1}}))_{A_{n}}$.
\end{proof} 
 
\begin{lem}\label{lem-dehn-free}
Let $M$ be a surface with $\kappa(M)\geq 0$ and let $\alpha, \alpha'\in V(C(M))$ with $i(\alpha, \alpha')\neq 0$. Then $t_{\alpha}^{n}$ and $t_{\alpha'}^{m}$ generate a free group of rank $2$ for all sufficiently large $n, m\in \mathbb{N}$, where $t_{\alpha}, t_{\alpha'}\in \Gamma(M)$ denote the Dehn twists about $\alpha$, $\alpha'$, respectively.
\end{lem}

\begin{proof}
We regard $\alpha$ and $\alpha'$ as elements in $\mathcal{PMF}$. Choose an open neighborhood $U$ of $\alpha$ such that
\[\overline{U}\subset \{ F\in \mathcal{PMF}: i(F, \alpha')\neq 0\},\]
where $\overline{K}$ denotes the closure of a subset $K$ of $\mathcal{PMF}$. Choose an open neighborhood $U'$ of $\alpha'$ such that
\[\overline{U'}\subset \{ F\in \mathcal{PMF}: i(F, \alpha)\neq 0\}\]
and $\overline{U}\cap \overline{U'}=\emptyset$. It follows from \cite[Theorem 4.3]{ivanov1} that there exist $n, m\in \mathbb{N}$ such that
\[t_{\alpha}^{k}(\overline{U'})\subset U\subset \overline{U},\ \ t_{\alpha'}^{l}(\overline{U})\subset U'\subset \overline{U'}\]
for any $k, l\in \mathbb{Z}$ with $|k|\geq n$, $|l|\geq m$. 

The lemma follows from the above inclusions and the following ping-pong argument: we show that $a=t_{\alpha}^{n}$ and $b=t_{\alpha'}^{m}$ generate a free group of rank $2$. Let $w$ be a non-empty reduced word consisting of $a^{\pm 1}$, $b^{\pm 1}$. We prove that $w$ is non-trivial in $\Gamma(M)$. It follows from the above inclusion that both $a^{k}$ and $b^{l}$ are non-trivial for any $k, l\in \mathbb{Z}\setminus \{ 0\}$. Therefore, by possibly replacing $w$ by an appropriate conjugate and an inverse, it is enough to prove that $w=a^{k}w'b^{l}$ is non-trivial in $\Gamma(M)$, where $k, l\in \mathbb{Z}\setminus \{ 0\}$ and $w'$ is a reduced word such that if $w'$ is non-empty, then the first letter of $w'$ is $b$ or $b^{-1}$ and the last letter of $w'$ is $a$ or $a^{-1}$. Then $w(x)\in U$ for any $x\in \overline{U}\setminus U$ by the above inclusion, and in particular, $w(x)\neq x$, which implies that $w$ is non-trivial in $\Gamma(M)$.  
\end{proof}

\begin{lem}\label{lem-Psi2}
If $\alpha, \alpha'\in V(C(M_{1}))$ satisfy $i(\alpha, \alpha')\neq 0$, then $i(\Psi(x, \alpha), \Psi(x, \alpha'))\neq 0$ for a.e.\ $x\in Y_{1}$.
\end{lem}

\begin{proof}
Let $\alpha, \alpha'\in V(C(M_{1}))$ with $i(\alpha, \alpha')\neq 0$. Assume that there exists a Borel subset $A$ of $Y_{1}$ with positive measure satisfying the following conditions: 
\begin{enumerate}
\item[(i)] $\Psi(\cdot, \alpha)$ and $\Psi(\cdot, \alpha')$ are constant on $A$. Let $\beta, \beta'\in V(C(M_{2}))$ be their values, respectively;
\item[(ii)] $i(\beta, \beta')=0$ and 
\[f((\mathcal{G}_{\alpha}^{1})_{A})<(\mathcal{G}_{\beta}^{2})_{f(A)},\ \ f((\mathcal{G}_{\alpha'}^{1})_{A})<(\mathcal{G}_{\beta'}^{2})_{f(A)}.\]
\end{enumerate}
Since $i(\beta, \beta')=0$, we see that $(\mathcal{G}_{\beta}^{2})_{f(A)}\vee (\mathcal{G}_{\beta'}^{2})_{f(A)}$ is amenable. On the other hand, $(\mathcal{G}_{\alpha}^{1})_{A}\vee (\mathcal{G}_{\alpha'}^{1})_{A}$ is non-amenable by Lemma \ref{lem-dehn-free} and Lemma \ref{lem-free-prod}, which is a contradiction. 
\end{proof} 
 
For each $\alpha \in V(C(M_{1}))$, we have a Borel subset $A_{\alpha}$ of $Y_{1}$ with full measure such that $\Psi(\cdot, \alpha)$ is defined on $A_{\alpha}$. Put \[A_{1}=\bigcap_{\alpha \in V(C(M_{1}))}A_{\alpha}.\]
By Lemmas \ref{lem-Psi1} and \ref{lem-Psi2}, for each pair $\{ \alpha, \alpha'\}$ of elements in $V(C(M_{1}))$, we can take a Borel subset $A_{\alpha, \alpha'}$ of $A_{1}$ with full measure so that for any $x\in A_{\alpha, \alpha'}$, we have $i(\Psi(x, \alpha), \Psi(x, \alpha'))=0$ if $i(\alpha, \alpha')=0$ and $i(\Psi(x, \alpha), \Psi(x, \alpha'))\neq 0$ if $i(\alpha, \alpha')\neq 0$. Put
\[A=\bigcap_{\alpha, \alpha'\in V(C(M_{1}))}A_{\alpha, \alpha'}.\]
Then $\Psi(x, \alpha)$ is defined for any $x\in A$ and $\alpha \in V(C)$, and both Lemmas \ref{lem-Psi1} and \ref{lem-Psi2} are satisfied for any $x\in A$.

 

Under the assumption $(\bullet)$, suppose that the two surfaces $M_{1}$, $M_{2}$ are equal and we denote the surface by $M$. Applying the above process to $f$ and $f^{-1}$, we see that there exist a Borel subset $A$ of $Y_{1}$ with full measure and a Borel map
\[\Psi \colon A\times V(C)\rightarrow V(C)\]
such that for each $x\in A$, the map $\Psi(x, \cdot)\colon V(C)\rightarrow V(C)$ defines an element of the automorphism group ${\rm Aut}(C)$ of the curve complex. We can define a Borel map $\Psi \colon A\rightarrow {\rm Aut}(C)$ by $\Psi(x)=\Psi(x, \cdot)$ for $x\in A$.  
 
For simplicity, we denote $\pi \circ \rho_{i}$ by $\rho_{i}$ for $i=1, 2$, where $\pi \colon \Gamma \rightarrow {\rm Aut}(C)$ is the natural homomorphism. 
 
\begin{lem}\label{lem-main-eq}
We have the equality
\[\Psi(r(\gamma))=\rho_{2}(f(\gamma))\Psi(s(\gamma))\rho_{1}(\gamma^{-1})\]
for a.e.\ $\gamma \in (\mathcal{G}^{1})_{Y_{1}}$.
\end{lem}




\begin{proof}
Let $A$ be a Borel subset of $Y_{1}$ and let $g_{1}\in \Gamma_{1}$, $g_{2}\in \Gamma_{2}$ be elements satisfying the following conditions: 
\begin{enumerate}
\item[(a)] $(g_{1}, x)\in (\mathcal{G}^{1})_{Y_{1}}$ and $(g_{2}, f(x))=f(g_{1}, x)\in (\mathcal{G}^{2})_{Y_{2}}$ for any $x\in A$;
\item[(b)] the map $\Psi$ is constant on $A$ and $g_{1}A$, respectively. Let $\psi, \psi'\in {\rm Aut}(C)$ be the values on $A$ and $g_{1}A$, respectively. 
\end{enumerate}
For each $\alpha \in V(C)$, there exists a Borel subset $B$ of $A$ with positive measure such that $f((\mathcal{G}_{\alpha}^{1})_{B})<(\mathcal{G}_{\psi(\alpha)}^{2})_{f(B)}$. Note that for $\alpha \in V(C)$ and $g\in \Gamma(M)$, we have
\[gt_{\alpha}g^{-1}=t_{g\alpha}\]
by \cite[Lemma 4.1.C]{ivanov2}, where $t_{\beta}\in \Gamma(M)$ denotes the Dehn twist about $\beta \in V(C)$. It follows that 
\[(g_{1}, r(\gamma))\gamma (g_{1}^{-1}, g_{1}s(\gamma))\in (\mathcal{G}_{g_{1}\alpha}^{1})_{g_{1}B},\ \ (g_{2}, r(\delta))\delta (g_{2}^{-1}, g_{2}s(\delta))\in (\mathcal{G}_{g_{2}\psi(\alpha)}^{2})_{f(g_{1}B)}\]
for $\gamma \in (\mathcal{G}_{\alpha}^{1})_{B}$ and $\delta \in (\mathcal{G}_{\psi(\alpha)})_{f(B)}$. Therefore, $f((\mathcal{G}_{g_{1}\alpha}^{1})_{g_{1}B})<(\mathcal{G}_{g_{2}\psi(\alpha)}^{2})_{f(g_{1}B)}$. Thus, $\psi'(g_{1}\alpha)=g_{2}\psi(\alpha)$. Since this equality holds for any $\alpha \in V(C)$, we have $\psi'=g_{2}\psi g_{1}^{-1}$. This means that
\[\Psi(r(\gamma))=\rho_{2}(f(\gamma))\Psi(s(\gamma))\rho_{1}(\gamma)^{-1}\]
for a.e.\ $\gamma =(g_{1}, x)\in (\mathcal{G}^{1})_{Y_{1}}$ with $x\in A$. 
\end{proof}

\begin{defn}
Let $S$ be a Borel space and let $m$ be a measure on $S$.
\begin{enumerate}
\item[(i)] Suppose that we are given a Borel space $T$, a Borel action of a discrete group $G$ on $S$, $T$ and a Borel map $f\colon S\rightarrow T$. We say that the map $f$ is almost $G$-equivariant if the equality
\[f(gx)=gf(x)\]
holds for any $g\in G$ and a.e.\ $x\in S$. 
\item[(ii)] Suppose that we have discrete groups $\Gamma$, $\Lambda$, $G$ and homomorphisms $\pi \colon \Gamma \rightarrow G$, $\tau \colon \Lambda \rightarrow G$. Then we denote by $(G, \pi, \tau)$ the Borel space $G$ equipped with the $\Gamma \times \Lambda$-action given by
\[(\gamma, \lambda)g=\pi(\gamma)g\tau(\lambda)^{-1}\]
for $\gamma \in \Gamma$, $\lambda \in \Lambda$ and $g\in G$.
\end{enumerate}
\end{defn}

\begin{thm}
For $i=1, 2$, let $\Gamma_{i}$ be a finite index subgroup of $\Gamma(M;m_{i})$, where $M$ is a surface with $\kappa(M)>0$ and $m_{i}\geq 3$ is an integer. Suppose that we have a ME coupling $(\Omega, m)$ of $\Gamma_{1}$ and $\Gamma_{2}$. Then there exists an essentially unique almost $\Gamma_{1}\times \Gamma_{2}$-equivariant Borel map $\Phi \colon \Omega \rightarrow ({\rm Aut}(C), \pi, \pi)$, where $\pi \colon \Gamma(M)^{\diamond}\rightarrow {\rm Aut}(C)$ is the natural homomorphism.
\end{thm}

\begin{proof}
As in Subsection \ref{subsec-ME oe}, we can construct a measure-preserving $\Gamma_{i}$-action on a standard finite measure space $(X_{i}, \mu_{i})$ for $i=1, 2$ such that they satisfy the assumption $(\bullet)$. In this proof, we use the notation as in $(\bullet)$. For the existence of $\Phi$, it is enough to show that there exists a Borel map $\Phi \colon \Omega \rightarrow {\rm Aut}(C)$ such that
\[\Phi(g_{1}g_{2}\omega)=\pi(g_{2})\Phi(\omega)\pi(g_{1})^{-1}\]
for any $g_{1}\in \Gamma_{1}$, $g_{2}\in \Gamma_{2}$ and a.e.\ $\omega \in \Omega$. By Lemma \ref{lem-ME iso}, the space $\Omega$ is isomorphic to $X_{1}\times \Gamma_{2}$ as $\Gamma_{1}\times \Gamma_{2}$-spaces. Here, the $\Gamma_{1}\times \Gamma_{2}$-action on $X_{1}\times \Gamma_{2}$ is given by the formula
\[g_{1}g_{2}(x, \gamma)=(g_{1}x, \alpha(g_{1}, x)\gamma g_{2}^{-1})\]
for $g_{1}\in \Gamma_{1}$, $g_{2}, \gamma \in \Gamma_{2}$ and $x\in X_{1}$, where $\alpha \colon \Gamma_{1}\times X_{1}\rightarrow \Gamma_{2}$ is the associated cocycle. We identify $\Omega$ with $X_{1}\times \Gamma_{2}$. For the proof of the theorem, it is enough to show that if we define $\Phi \colon \Omega \rightarrow {\rm Aut}(C)$ by the formula
\[\Phi(g_{1}g_{2}(x, e))=\pi(g_{2})\Psi(x)\pi(g_{1})^{-1}\]
for $g_{1}\in \Gamma_{1}$, $g_{2}\in \Gamma_{2}$ and $x\in Y_{1}$, then it is well-defined. In other words, it is enough to show that 
\[\pi(g_{2})\Psi(x)\pi(g_{1})^{-1}=\pi(g_{2}')\Psi(x')\pi(g_{1}')^{-1}\]
for any $g_{1}, g_{1}'\in \Gamma_{1}$, $g_{2}, g_{2}'\in \Gamma_{2}$ and a.e.\ $x, x'\in Y_{1}$ satisfying 
\[g_{1}g_{2}(x, e)=g_{1}'g_{2}'(x', e).\]
In what follows, we omit $\pi$ for simplicity. Since 
\[(x', e)=(g_{1}')^{-1}g_{1}(g_{2}')^{-1}g_{2}(x, e)=((g_{1}')^{-1}g_{1}x, \alpha((g_{1}')^{-1}g_{1}, x)g_{2}^{-1}g_{2}'),\]
we see that $x'=(g_{1}')^{-1}g_{1}x\in Y_{1}$. Since $\alpha(g, y)=\rho_{2}(f(g, y))$ for $g\in \Gamma_{1}$, $y\in Y_{1}$ with $gy\in Y_{1}$, we have
\begin{align*}
\Psi(x')&=\Psi((g_{1}')^{-1}g_{1}x)=\rho_{2}(f((g_{1}')^{-1}g_{1}, x))\Psi(x)\rho_{1}((g_{1}')^{-1}g_{1}, x)^{-1}\\
&=(g_{2}')^{-1}g_{2}\Psi(x)g_{1}^{-1}g_{1}'
\end{align*}
by Lemma \ref{lem-main-eq}, which shows the claim. 

The uniqueness of $\Phi$ is a consequence of Theorem \ref{thm-icc} and the following Lemma \ref{lem-gen-uni}.
\end{proof}
 
\begin{defn}
Let $\pi \colon \Gamma \rightarrow G$ be a homomorphism between discrete groups. Then $\pi$ is said to be ICC (infinite conjugacy class) if the set $\{ \pi(\gamma)g\pi(\gamma)^{-1}: \gamma \in \Gamma \}$
consists of infinitely many elements for any $g\in G\setminus \{ e\}$.
\end{defn}

\begin{lem}\label{lem-gen-uni}
Let $\Gamma$, $\Lambda$ and $G$ be discrete groups and assume that
\[\pi \colon \Gamma \rightarrow G,\ \ \tau \colon \Lambda \rightarrow G\]
are homomorphisms and that either $\pi$ or $\tau$ is ICC. Suppose the following two conditions:
\begin{enumerate}
\item[(i)] we have a ME coupling $(\Sigma, m)$ of $\Gamma$ and $\Lambda$;
\item[(ii)] there exist two almost $\Gamma \times \Lambda$-equivariant Borel maps $\Phi, \Phi' \colon \Sigma \rightarrow (G, \pi, \tau)$.
\end{enumerate}
Then $\Phi$ and $\Phi'$ are essentially equal. 
\end{lem}

\begin{proof}
We may assume that $\pi$ is ICC. Let us define a Borel map $\Phi_{0}\colon \Sigma \rightarrow G$ by $\Phi_{0}(x)=\Phi'(x)\Phi(x)^{-1}$ for $x \in \Sigma$. Then $\Phi_{0}$ satisfies the equality
\[\Phi_{0}(\gamma \lambda x)=\pi(\gamma)\Phi_{0}(x)\pi(\gamma)^{-1}\]
for any $\gamma \in \Gamma$, $\lambda \in \Lambda$ and a.e.\ $x \in \Sigma$. Therefore, $\Phi_{0}$ is $\Lambda$-invariant and induces an almost $\Gamma$-equivariant Borel map $\Lambda \backslash \Sigma \rightarrow G$, where the $\Gamma$-action on $G$ is given by conjugation via $\pi$. By projecting the finite $\Gamma$-invariant measure on $\Lambda \backslash \Sigma$ to $G$, we obtain a finite measure on $G$ which is invariant under the conjugation via $\pi$ of each element of $\Gamma$. Since $\pi$ is ICC, the support of this measure is equal to $\{ e\}$, which means that $\Phi_{0}=e$ a.e.\ on $\Sigma$. 
\end{proof}

\begin{lem}\label{lem-gen-ME eq-ext}
Let $\Gamma$, $\Lambda$ and $G$ be discrete groups and assume that
\[\pi \colon \Gamma \rightarrow G,\ \ \tau \colon \Lambda \rightarrow G\]
are homomorphisms. Suppose the following three conditions:
\begin{enumerate}
\item[(i)] we have a normal subgroup $\Gamma'$ of $\Gamma$ (resp. $\Lambda'$ of $\Lambda$) of finite index and a ME coupling $(\Sigma, m)$ of $\Gamma$ and $\Lambda$;
\item[(ii)] either the restrictions $\pi \colon \Gamma'\rightarrow G$ or $\tau \colon \Lambda'\rightarrow G$ is ICC;
\item[(iii)] there exists an almost $\Gamma' \times \Lambda'$-equivariant Borel map $\Phi \colon \Sigma \rightarrow (G, \pi, \tau)$.
\end{enumerate}
Then the map $\Phi$ is almost $\Gamma \times \Lambda$-equivariant. 
\end{lem}
 
\begin{proof}
We may assume that the restriction $\tau \colon \Lambda'\rightarrow G$ is ICC. For fixed $\gamma \in \Gamma$ and $\lambda \in \Lambda$, define a Borel map $\Phi_{0}\colon \Sigma \rightarrow G$ by the formula 
\[\Phi_{0}(x)=\Phi(\gamma \lambda x)^{-1}\pi(\gamma) \Phi(x)\tau(\lambda)^{-1}\]
for $x\in \Sigma$. 

Let $g \in \Gamma'$, $h\in \Lambda'$. Since $\Gamma'$ is normal in $\Gamma$ and $\Lambda'$ is normal in $\Lambda$, we have $g'\in \Gamma'$ and $h'\in \Lambda'$ such that $\gamma g=g'\gamma$ and $\lambda h'=h\lambda$. Then
\begin{align*}
\Phi_{0}(gh'x)&=\Phi(\gamma g\lambda h'x)^{-1}\pi(\gamma)\pi(g)\Phi(x)\tau(h')^{-1}\tau(\lambda)^{-1}\\
                  &=\Phi(g'\gamma h\lambda x)^{-1}\pi(g')\pi(\gamma)\Phi(x)\tau(\lambda)^{-1}\tau(h)^{-1}\\
                  &=\tau(h)\Phi(\gamma \lambda x)^{-1}\pi(\gamma)\Phi(x)\tau(\lambda)^{-1}\tau(h)^{-1}\\
                  &=\tau(h)\Phi_{0}(x)\tau(h)^{-1}.
\end{align*}
Since $g\in \Gamma'$ is arbitrary, the map $\Phi_{0}$ induces a Borel map $\Gamma'\backslash \Sigma \rightarrow G$. The projectivized finite measure on $G$ is invariant under the conjugation via $\tau$ of each element of $\Lambda'$. As in the proof of Lemma \ref{lem-gen-uni}, we can show that $\Phi_{0}=e$ a.e.\ on $\Sigma$.
\end{proof}

\begin{cor}\label{cor-mcg-ME eq}
Let $M$ be a surface with $\kappa(M)>0$ and let $\Gamma$, $\Lambda$ be finite index subgroups of $\Gamma(M)^{\diamond}$ or ${\rm Aut}(C)$, respectively. Suppose that we have a ME coupling $(\Omega, m)$ of $\Gamma$ and $\Lambda$. Then there exists an essentially unique almost $\Gamma \times \Lambda$-equivariant Borel map $\Omega \rightarrow ({\rm Aut}(C), \pi, \pi)$. 
\end{cor}

 
\section{Measure equivalence rigidity}

Given a ME coupling $(\Sigma, m)$ of discrete groups $\Gamma$ and $\Lambda$, we can define the opposite coupling $\check{\Sigma}$ of $\Lambda$ and $\Gamma$ as the $\Lambda \times \Gamma$-space obtained by the canonical isomorphism between $\Gamma \times \Lambda$ and $\Lambda \times \Gamma$. When $\Gamma =\Lambda$, we distinguish the $\Gamma$-actions on $\Sigma$ by writing $(\gamma, x)\mapsto A_{\gamma}^{1}x$ and $A_{\gamma}^{2}x$, respectively, for $\gamma \in \Gamma$ and $x\in \Sigma$. 

If $(\Sigma, m)$ is a ME coupling of discrete groups $\Gamma$ and $\Lambda$ and $(\Omega, n)$ is a ME coupling of discrete groups $\Lambda$ and $\Delta$, then we can define the composed coupling $\Sigma \times_{\Lambda}\Omega$ of $\Gamma$ and $\Delta$ as the quotient space of $\Sigma \times \Omega$ by the diagonal $\Lambda$-action.

\begin{defn}
Let $\pi \colon \Gamma \rightarrow G$ be a homomorphism between discrete groups. We say that $\pi$ is almost an isomorphism if both $\ker(\pi)$ and ${\rm coker}(\pi)$ are finite.
\end{defn}

\begin{thm}\label{thm-gen-ME eq}
Let $\Gamma$, $\Lambda$, $G$ be discrete groups and let $\pi, \tau \colon \Gamma \rightarrow G$ be homomorphisms. Suppose that $\pi$ is ICC and $\tau$ is almost an isomorphism and we have a ME coupling $(\Sigma, m)$ of $\Gamma$ and $\Lambda$ and let $\Omega =\Sigma \times_{\Lambda}\Lambda \times_{\Lambda}\check{\Sigma}$ be the self ME coupling of $\Gamma$. Moreover, assume that there exists an almost $\Gamma \times \Gamma$-equivariant Borel map $\Phi \colon \Omega \rightarrow (G, \pi, \tau)$. Then we can find the following two maps:
\begin{enumerate}
\item[(a)] a homomorphism $\rho \colon \Lambda \rightarrow G$ which is almost an isomorphism;
\item[(b)] an almost $\Gamma \times \Lambda$-equivariant Borel map $\Phi_{0} \colon \Sigma \rightarrow (G, \tau, \rho)$. 
\end{enumerate}
\end{thm}

Before the proof, we give the following lemma:

\begin{lem}\label{lem-cocycle-finite}
Let $\Gamma$, $\Lambda$, $G$ be discrete groups and let $(\Sigma, m)$ be a ME coupling of $\Gamma$ and $\Lambda$. Let $Y\subset \Sigma$ be a fundamental domain of the $\Gamma$-action on $\Sigma$ and let $\theta \colon \Lambda \times Y\rightarrow \Gamma$ be the associated cocycle. Suppose the following two conditions:
\begin{enumerate}
\item[(i)] we have a homomorphism $\pi \colon \Gamma \rightarrow G$ which is almost an isomorphism;
\item[(ii)] there exists a subgroup $G_{0}$ of $G$ such that the cocycle $\pi \circ \theta \colon \Lambda \times Y\rightarrow G$ is cohomologous to a cocycle which is essentially valued in $G_{0}$.  
\end{enumerate}
Then $G_{0}$ is a subgroup of finite index in $G$.
\end{lem}

\begin{proof}
Take a standard $\Lambda$-action on a standard probability space $X_{0}$. Define a $\Gamma$-action and a $G\times \Lambda$-action on $\Sigma \times G\times X_{0}$ by 
\[\gamma(z, g, x)=(\gamma z, \pi(\gamma)g, x)\]
\[(g_{1}, \lambda)(z, g, x)=(\lambda z, gg_{1}^{-1}, \lambda x)\]
for $g, g_{1}\in G$, $\gamma \in \Gamma$, $\lambda \in \Lambda$, $z\in \Sigma$ and $x\in X_{0}$. Consider the quotient $G\times \Lambda$-space $\tilde{\Sigma}$ of $\Sigma \times G\times X_{0}$ by the $\Gamma$-action. Since $\ker(\pi)$ is finite, the $\Lambda$-action has a fundamental domain. Since ${\rm coker}(\pi)$ is finite, the $\Lambda$-action has a fundamental domain of finite measure. Thus, the $G\times \Lambda$-space $\tilde{\Sigma}$ is a ME coupling of $G$ and $\Lambda$. 

Let $p\colon \Sigma \times G\times X_{0}\rightarrow \tilde{\Sigma}$ be the natural projection. Then $p(Y\times \{ e\} \times X_{0})$ is a fundamental domain of the $G$-action on $\tilde{\Sigma}$. Remark that $p$ is injective on $Y\times \{ e\} \times X_{0}$. The cocycle $\tilde{\theta}\colon \Lambda \times p(Y\times \{ e\} \times X_{0})\rightarrow G$ associated to it is given by 
\[\tilde{\theta}(\lambda, p(y, e, x))=\pi\circ \theta(\lambda, y)\]
for $\lambda \in \Lambda$, $y\in Y$ and $x\in X_{0}$. By assumption, we can find a Borel map $\varphi \colon Y\rightarrow G$ such that 
\[\theta'(\lambda, y)=\varphi(\lambda \cdot y)\pi \circ \theta(\lambda, y)\varphi(y)^{-1}\in G_{0}\]
for any $\lambda \in \Lambda$ and a.e.\ $y\in Y$. Define a Borel map $\tilde{\varphi}\colon p(Y\times \{ e\} \times X_{0})\rightarrow G$ by $\tilde{\varphi}(p(y, e, x))=\varphi(y)$ for $y\in Y$. Then
\begin{align*}
&\tilde{\varphi}(\lambda \cdot p(y, e, x))\tilde{\theta}(\lambda, p(y, e, x))\tilde{\varphi}(p(y, e, x))^{-1}\\
=\ &\varphi(\lambda \cdot y)\pi \circ \theta(\lambda, y)\varphi(y)^{-1}\in G_{0},
\end{align*}
which means that $\tilde{\theta}$ is cohomologous to a cocycle which is essentially valued in $G_{0}$. The lemma now follows from \cite[Lemma 6.1]{ms}.  
\end{proof}

\vspace{1em}

\noindent {\it Proof of Theorem \ref{thm-gen-ME eq}.}
This proof is almost the same as in \cite[Section 6.2]{ms}. One denotes the element corresponding to $(x, \lambda, y)\in \Sigma \times \Lambda \times \check{\Sigma}$ by $[x, \lambda, y]\in \Sigma \times_{\Lambda}\Lambda \times_{\Lambda}\check{\Sigma}$. As in \cite[Lemma 6.6]{ms}, we can prove the following lemma by using the assumption that $\pi$ is ICC:

\begin{lem}
If one defines a Borel map $\Psi \colon \Sigma^{3}\rightarrow G$ by
\[\Psi(x, y, z)=\Phi([x, e, z])\Phi([y, e, z])^{-1}\]
for $(x, y, z)\in \Sigma^{3}$, then 
\[\Psi(x, y, z_{1})=\Psi(x, y, z_{2})\]
for $m^{4}$-a.e.\ $(x, y, z_{1}, z_{2})\in \Sigma^{4}$.
\end{lem} 

Define a Borel map $F\colon \Sigma^{2}\rightarrow G$ by $F(x, y)=\Psi(x, y, z)$. It follows from \cite[Lemma 6.2]{ms} that for $m$-a.e.\ $x\in \Sigma$, the Borel map $\rho_{x}\colon \Lambda \rightarrow \Gamma$ given by
\[\rho_{x}(\lambda)=F(\lambda^{-1}x, y)F(x, y)^{-1}\]
is the same for $m$-a.e.\ $y\in \Sigma$ and defines a homomorphism. Moreover, the equality 
\[\rho_{y}(\lambda)=F(x, y)^{-1}\rho_{x}(\lambda)F(x, y)\]
holds for any $\lambda \in \Lambda$ and $m^{2}$-a.e.\ $(x, y)\in \Sigma^{2}$. Note that we have the equality
\[\rho_{x}(\lambda)=\Phi([x, \lambda, z])\Phi([x, e, z])^{-1}\]
for any $\lambda \in \Lambda$ and $m^{2}$-a.e.\ $(x, z)\in \Sigma^{2}$. Let $N$ be the normal subgroup of $\Lambda$ which is the common kernel of $\rho_{x}$ for $m$-a.e.\ $x\in \Sigma$.

Let $D\subset \Sigma$ be a fundamental domain of the $\Lambda$-action on $\Sigma$ and put $\tilde{\Omega}=D\times \Lambda \times D\subset \Sigma \times \Lambda \times \Sigma$. This inclusion induces a Borel  isomorphism between $\tilde{\Omega}$ and $\Omega$. Let us give a $\Gamma$-action on $\tilde{\Omega}$ induced by the second $\Gamma$-action $A^{2}$ on $\Omega$ and give a $\Lambda$-action on $\tilde{\Omega}$ by the left multiplication on the second coordinate:
\[(\gamma, \lambda)(x, \lambda_{1}, y)=(x, \lambda \lambda_{1}\alpha(\gamma, y)^{-1}, \gamma \cdot y)\]
for $\gamma \in \Gamma$, $\lambda, \lambda_{1}\in \Lambda$ and $x, y\in D$, where $\alpha \colon \Gamma \times D\rightarrow \Lambda$ is the associated cocycle to $D$. 

Let $\tilde{\Phi}\colon \tilde{\Omega}\rightarrow G$ be the Borel map induced by $\Phi$ and the isomorphism between $\tilde{\Omega}$ and $\Omega$. Note that $\tilde{\Phi}$ is almost $\Gamma$-equivariant in the following sense:
\[\tilde{\Phi}(\gamma \omega)=\tilde{\Phi}(\omega)\tau(\gamma)^{-1}\]
for any $\gamma \in \Gamma$ and a.e.\ $\omega \in \tilde{\Omega}$. Put $E_{0}=\tilde{\Phi}^{-1}(\{ g_{n}\} )$, where $\{ g_{n}\} \subset G$ is a finite set of all representatives of $G/\tau(\Gamma)$. Remark that $E_{0}$ is invariant under the action of $\ker(\tau)$. If $E\subset E_{0}$ is a fundamental domain of the $\ker(\tau)$-action on $E_{0}$, then it is also a fundamental domain of the $\Gamma$-action on $\tilde{\Omega}$. Since $\ker(\tau)$ is finite, the measure of $E_{0}$ is finite. The homomorphism $\rho_{x}$ is given by
\[\rho_{x}(\lambda)=\tilde{\Phi}(x, \lambda, z)\tilde{\Phi}(x, e, z)^{-1}\]
for any $\lambda \in \Lambda$ and $m^{2}$-a.e.\ $(x, z)\in D^{2}$.

For $\lambda_{0}\in \Lambda$, it is easy to see that $\lambda_{0}\in N$ if and only if
\[\tilde{\Phi}(x, \lambda_{0}\lambda_{1}, y)=\tilde{\Phi}(x, \lambda_{1}, y)\]
for any $\lambda_{1}\in \Lambda$ and $m^{2}$-a.e.\ $(x, y)\in D^{2}$. It follows that any element in $N$ preserves $E_{0}$. Since the measure of $E_{0}$ is finite, we see that $N$ is finite. Note that for any $\lambda \in \Lambda$ and a.e.\ $t=(x, \lambda_{1}, z)\in \tilde{\Omega}$, we have
\begin{align*}
\rho_{x}(\lambda)&=\rho_{x}(\lambda \lambda_{1})\rho_{x}(\lambda_{1})^{-1}\\
                         &=\tilde{\Phi}(x, \lambda \lambda_{1}, z)\tilde{\Phi}(x, e, z)^{-1}(\tilde{\Phi}(x, \lambda_{1}, z)\tilde{\Phi}(x, e, z)^{-1})^{-1}\\
                         &=\tilde{\Phi}(\lambda t)\tilde{\Phi}(t)^{-1}.
\end{align*}
Let $\theta \colon \Lambda \times E\rightarrow \Gamma$ be the associated cocycle to $E$. It follows from Fubini's theorem that there exists $x_{0}\in D$ such that $\rho=\rho_{x_{0}}\colon \Lambda \rightarrow G$ is a homomorphism with kernel $N$ and for any $\lambda \in \Lambda$, $\gamma \in \Gamma$, a.e.\ $x\in D$ and a.e.\ $(\lambda_{1}, z)\in \Lambda \times D$, we have 
\[\rho_{x}(\lambda)=F(x_{0}, x)^{-1}\rho(x)F(x_{0}, x)\]
and
\[\tilde{\Phi}((\gamma, \lambda)(x_{0}, \lambda_{1}, z))=\rho(\lambda)\tilde{\Phi}(x_{0}, \lambda_{1}, z)\tau(\gamma)^{-1}.\]

We show that $\rho(\Lambda)$ is a subgroup of finite index in $G$. Define a Borel map $\varphi \colon E\rightarrow G$ by
\[\varphi(t)=F(x_{0}, x)\tilde{\Phi}(t)\]
for $t=(x, \lambda_{1}, z)\in E$. Put
\[\theta'(\lambda, t)=\varphi(\lambda \cdot t)\tau \circ \theta(\lambda, t)\varphi(t)^{-1}\]
for $\lambda \in \Lambda$ and $t\in E$. Since $\lambda \cdot t=\lambda \theta(\lambda, t)t$, we see that
\[\tilde{\Phi}(\lambda \cdot t)=\tilde{\Phi}(\lambda \theta(\lambda, t)t)=\tilde{\Phi}(\lambda t)\tau \circ \theta(\lambda, t)^{-1}\]
and
\begin{align*}
\theta'(\lambda, t)&=F(x_{0}, x)\tilde{\Phi}(\lambda \cdot t)\tau \circ \theta(\lambda, t)\tilde{\Phi}(t)^{-1}F(x_{0}, x)^{-1}\\
&=F(x_{0}, x)\tilde{\Phi}(\lambda t)\tilde{\Phi}(t)^{-1}F(x_{0}, x)^{-1}\\
&=F(x_{0}, x)\rho_{x}(\lambda)F(x_{0}, x)^{-1}=\rho(\lambda)\in \rho(\Lambda).
\end{align*}
It follows from Lemma \ref{lem-cocycle-finite} that $\rho(\Lambda)$ is a subgroup of finite index in $G$.

Finally, we construct a Borel map $\Phi_{0}\colon \Sigma \rightarrow G$. Note that $\{ x_{0}\} \times \Lambda \times D\subset \tilde{\Omega}$ is a $\Gamma \times \Lambda$-invariant Borel subset isomorphic to $\Sigma$ as $\Gamma \times \Lambda$-spaces. It follows from the choice of $x_{0}$ that the composition of the restriction of $\tilde{\Phi}$ to $\{ x_{0}\} \times \Lambda \times D$ and the map $G\ni g\mapsto g^{-1}\in G$ is a desired map. \hfill $\Box$

\vspace{1em}

Combining Corollary \ref{cor-mcg-ME eq} and Theorem \ref{thm-gen-ME eq}, we obtain Theorem \ref{thm-main}.

\vspace{1em}

\noindent {\it Proof of Theorem \ref{cor-classification}.}
First, suppose that $\kappa(M^{1})\geq \kappa(M^{2})$. We may assume that $\kappa(M^{1})\geq 2$. By Theorem \ref{thm-main}, we can find an injective homomorphism $\Gamma(M^{1};3)\rightarrow {\rm Aut}(C(M^{2}))$ with finite cokernel. By using Theorem \ref{thm-cc-auto} and restricting the homomorphism to some subgroup $\Gamma_{1}$ of finite index in $\Gamma(M^{1};3)$, we can construct an injective homomorphism from $\Gamma_{1}$ into $\Gamma(M^{2})$ with finite cokernel. It follows from \cite[Theorem 2]{sha} that $M^{1}=M_{0, 6}$ and $M^{2}=M_{2, 0}$. Similarly, if we assume that $\kappa(M^{1})\leq \kappa(M^{2})$, then it can be shown that $M^{1}=M_{0, 6}$ and $M^{2}=M_{2, 0}$. \hfill $\Box$


\section{Rigidity for a direct product of mapping class groups}

We need to review Monod-Shalom's technique in \cite[Section 5.1]{ms}.

Let $\Gamma_{1}, \ldots, \Gamma_{n}$ be torsion-free discrete groups in the class $\mathcal{C}$, mentioned in Section \ref{sec-int}, and let $\Lambda_{1}, \ldots, \Lambda_{n}$ be torsion-free discrete groups. Put $\Gamma =\Gamma_{1}\times \cdots \times \Gamma_{n}$ and $\Lambda =\Lambda_{1}\times \cdots \times \Lambda_{n}$. Let us denote
\[\Gamma_{i}'=\prod_{j\neq i}\Gamma_{j},\ \ \Lambda_{i}'=\prod_{j\neq i}\Lambda_{j}\]
for $i\in \{ 1, \ldots, n\}$. Suppose that we have a ME coupling $(\Sigma, m)$ of $\Gamma$ and $\Lambda$. 

In the above situation, we can find a bijection $t\colon \{ 1, \ldots, n\} \rightarrow \{ 1, \ldots, n\}$ and fundamental domains $Y, X\subset \Sigma$ of the $\Gamma$-, $\Lambda$-actions on $\Sigma$, respectively, satisfying 
\[\Lambda_{t(i)}Y\subset \Gamma_{i}Y,\ \ \Gamma_{i}X\subset \Lambda_{t(i)}X\]
for any $i\in \{ 1, \ldots, n\}$. Let $\overline{\Sigma}_{i}$ be the space of ergodic components of the $\Gamma_{i}'\times \Lambda_{t(i)}'$-action on $(\Sigma, m)$ for $i\in \{ 1, \ldots, n\}$, which is naturally a $\Gamma_{i}\times \Lambda_{t(i)}$-space. Define a measure $\mu_{i}$ (resp. $\nu_{i}$) on $\overline{\Sigma}_{i}$ by projecting the restricted measure on $\Gamma_{i}Y$ (resp. $\Lambda_{t(i)}X$) through the natural map $\Gamma_{i}Y\rightarrow \overline{\Sigma}_{i}$ (resp. $\Lambda_{t(i)}X\rightarrow \overline{\Sigma}_{i}$). Then 
\begin{enumerate}
\item[(a)] $\mu_{i}$ and $\nu_{i}$ are absolutely continuous with respect to each other.
\item[(b)] both $\mu_{i}$ and $\nu_{i}$ are invariant for the $\Gamma_{i}\times \Lambda_{t(i)}$-action on $\overline{\Sigma}_{i}$.
\item[(c)] if $\overline{Y}$ (resp. $\overline{X}$) is the image of $Y$ (resp. $X$) in $\overline{\Sigma}_{i}$, then it is a fundamental domain of the $\Gamma_{i}$-action on $(\overline{\Sigma}_{i}, \mu_{i})$ (resp. the $\Lambda_{t(i)}$-action on $(\overline{\Sigma}_{i}, \nu_{i})$). Moreover, both $\mu_{i}(\overline{Y})$ and $\nu_{i}(\overline{X})$ are finite.
\end{enumerate}

These claims are shown in the proof of \cite[Theorem 1.16]{ms}, where the ergodicity of $(\Sigma, m)$ is assumed. However, we can show the above claims along the same line without this assumption.

Let 
\[c_{i}(x)=\frac{d\mu_{i}}{d\nu_{i}}(x),\ \ \ x\in \overline{\Sigma}_{i}\]
be the Radon-Nikodym derivative, which is positive and finite a.e.\ on $\overline{\Sigma}_{i}$. It follows from the condition (b) that the function $c_{i}$ is invariant for the $\Gamma_{i}\times \Lambda_{t(i)}$-action. Put
\[\overline{\Sigma}_{i, n}=\{ x\in \overline{\Sigma}_{i}: n<c_{i}(x)\leq n+1\}\]
for $n\in \mathbb{N}$. Then $\overline{\Sigma}_{i}=\bigsqcup_{n\in \mathbb{N}}\overline{\Sigma}_{i, n}$ up to null sets. It follows from the condition (c) that $\overline{\Sigma}_{i, n}$ is a ME coupling of $\Gamma_{i}$ and $\Lambda_{t(i)}$ with respect to $\mu_{i}$ for each $n\in \mathbb{N}$ (if $\overline{\Sigma}_{i, n}$ has non-zero measure). 

In this situation, we suppose the following condition: for $i\in \{ 1, \ldots, n\}$ and $j\in \{ 1, 2\}$, let $M_{i}^{j}$ be a surface with $\kappa(M_{i}^{j})>0$ and $M_{i}^{j}\neq M_{1, 2}, M_{2, 0}$. Assume that $\Gamma_{i}$ (resp. $\Lambda_{i}$) is a torsion-free subgroup of finite index in $\Gamma(M_{i}^{1})^{\diamond}$ (resp. $\Gamma(M_{i}^{2})^{\diamond}$) for each $i$.

Remark that the mapping class group $\Gamma(M)$ is in $\mathcal{C}$ for a surface $M$ with $\kappa(M)\geq 0$ (\cite[Corollary B]{ham2}) and this property is preserved under measure equivalence and in particular, commensurability up to finite kernel (\cite[Corollary 7.6]{ms}). Note that $M_{i}^{1}$ and $M_{t(i)}^{2}$ are diffeomorphic for any $i$ by Theorem \ref{cor-classification} and let $g_{i}$ be an isotopy class of a diffeomorphism $M_{t(i)}^{2}\rightarrow M_{i}^{1}$. Let 
\[\pi_{g}\colon \prod_{i=1}^{n}\Gamma(M_{i}^{2})^{\diamond}\rightarrow \prod_{i=1}^{n}{\rm Aut}(C(M_{i}^{1}))\]
be the isomorphism defined by
\[\pi_{g}(\gamma_{1}, \ldots, \gamma_{n})=(\pi(g_{1}\gamma_{t(1)}g_{1}^{-1}), \ldots, \pi(g_{n}\gamma_{t(n)}g_{n}^{-1}))\]
for $\gamma_{i}\in \Gamma(M_{i}^{2})^{\diamond}$, where we denote by the same symbol $\pi$ the natural homomorphism $\Gamma(M)^{\diamond}\rightarrow {\rm Aut}(C(M))$ for a surface $M$. By applying Corollary \ref{cor-mcg-ME eq} to each ME coupling $\overline{\Sigma}_{i, n}$ of $\Gamma_{i}$ and $\Lambda_{t(i)}$, we can find an almost $\Gamma_{i}\times \Lambda_{t(i)}$-equivariant Borel map $\Phi_{i}\colon \overline{\Sigma}_{i}\rightarrow ({\rm Aut}(C(M_{i})), \pi, \pi_{g_{i}})$, where $\pi_{g_{i}}\colon \Gamma(M^{2}_{t(i)})^{\diamond}\rightarrow {\rm Aut}(C(M_{i}^{1}))$ is defined by using $g_{i}$. Then we can define a Borel map $\Phi \colon \Sigma \rightarrow \prod_{i=1}^{n}{\rm Aut}(C(M_{i}^{1}))$ by
\[\Phi(x)=(\Phi_{1}(p_{1}(x)), \ldots, \Phi_{n}(p_{n}(x)))\]
for $x\in \Sigma$, where $p_{i}\colon \Sigma \rightarrow \overline{\Sigma}_{i}$ denotes the natural projection. It is easy to see that
\[\Phi((\gamma, \lambda)x)=\pi(\gamma)\Phi(x)\pi_{g}(\lambda)^{-1}\]
for any $\gamma \in \Gamma$, $\lambda \in \Lambda$ and a.e.\ $x\in \Sigma$.

Hence, we have shown the following:

\begin{thm}\label{thm-prod-ME eq}
For $i\in \{ 1, \ldots, n\}$ and $j\in \{ 1, 2\}$, let $M_{i}^{j}$ be a surface with $\kappa(M_{i}^{j})>0$ and $M_{i}^{j}\neq M_{1, 2}, M_{2, 0}$. Assume that $\Gamma_{i}$ (resp. $\Lambda_{i}$) is a torsion-free subgroup of finite index in $\Gamma(M_{i}^{1})^{\diamond}$ (resp. $\Gamma(M_{i}^{2})^{\diamond}$). Put $\Gamma =\Gamma_{1}\times \cdots \times \Gamma_{n}$, $\Lambda =\Lambda_{1}\times \cdots \times \Lambda_{n}$. Suppose that we have a ME coupling $(\Sigma, m)$ of $\Gamma$ and $\Lambda$. Then we can find the following:
\begin{enumerate}
\item[(a)] a bijection $t$ on the set $\{ 1,\ldots, n\}$;
\item[(b)] an isotopy class $g_{i}$ of a diffeomorphism $M_{t(i)}^{2}\rightarrow M_{i}^{1}$ for each $i$;
\item[(c)] an almost $\Gamma \times \Lambda$-equivariant Borel map 
\[\Phi \colon \Sigma \rightarrow \left(\prod_{i=1}^{n}{\rm Aut}(C(M_{i}^{1})), \pi, \pi_{g}\right).\]
\end{enumerate}
\end{thm}



\begin{cor}\label{cor-prod-me-eq}
The conclusion of Theorem \ref{thm-prod-ME eq} holds even if $\Gamma$ (resp. $\Lambda$) is a subgroup of finite index in $\Gamma(M_{1}^{1})^{\diamond}\times \cdots \times \Gamma(M_{n}^{1})^{\diamond}$ (resp. $\Gamma(M_{1}^{2})^{\diamond}\times \cdots \times \Gamma(M_{n}^{2})^{\diamond}$).
\end{cor}


\begin{proof}
It is easy to check that if $\kappa(M_{i})>0$ and $\Gamma_{i}$ is a finite index subgroup of $\Gamma(M_{i})^{\diamond}$, then the natural homomorphism 
\[\Gamma_{1}\times \cdots \times \Gamma_{n}\rightarrow {\rm Aut}(C(M_{1}))\times \cdots \times {\rm Aut}(C(M_{n}))\]
is almost an isomorphism and ICC. By Lemma \ref{lem-gen-ME eq-ext} and Theorem \ref{thm-prod-ME eq}, we obtain the corollary.  
\end{proof}

Combining Theorem \ref{thm-gen-ME eq} and the above corollary, we can obtain Theorem \ref{thm-main-prod}. The following corollary determines all isomorphisms between finite index subgroups of a direct product of mapping class groups (see also \cite[Section 8.5]{ivanov2}).

\begin{cor}
For $i\in \{ 1, \ldots, n\}$, let $M_{i}$ be a surface with $\kappa(M_{i})>0$ and $M_{i}\neq M_{1, 2}, M_{2, 0}$, and let $\Gamma$ be a finite index subgroup of $G=\Gamma(M_{1})^{\diamond}\times \cdots \times \Gamma(M_{n})^{\diamond}$. Suppose that we have an injective homomorphism $\tau \colon \Gamma \rightarrow G$ with finite cokernel. Then we can find a bijection $t$ on the set $\{ 1,\ldots, n\}$ and an isotopy class $g_{i}$ of a diffeomorphism $M_{t(i)}\rightarrow M_{i}$ for each $i$ such that for any $\gamma =(\gamma_{1}, \ldots, \gamma_{n})\in \Gamma$, we have
\[\tau(\gamma)=(g_{1}\gamma_{t(1)}g_{1}^{-1}, \ldots, g_{n}\gamma_{t(n)}g_{n}^{-1}).\]
\end{cor}

\begin{proof}
We identify $\Gamma(M_{i})^{\diamond}$ and ${\rm Aut}(C(M_{i}))$ via the natural isomorphism. Consider the self ME coupling $(G, \pi, \tau)$ of $\Gamma$. It follows from Corollary \ref{cor-prod-me-eq} that we can find the following:
\begin{enumerate}
\item[(a)] a bijection $t$ on $\{ 1,\ldots, n\}$;
\item[(b)] an isotopy class $g_{i}$ of a diffeomorphism $M_{t(i)}^{2}\rightarrow M_{i}^{1}$ for each $i$;
\item[(c)] an almost $\Gamma \times \Gamma$-equivariant Borel map 
\[\Phi \colon (G, \pi, \tau) \rightarrow (G, \pi, \pi_{g}).\]
\end{enumerate}
Put $h=(h_{1},\ldots , h_{n})=\Phi(e)$ and define an automorphism $\pi_{hg}$ of $G$ by
\[\pi_{hg}(s)=(h_{1}g_{1}s_{t(1)}g_{1}^{-1}h_{1}^{-1}, \ldots, h_{n}g_{n}s_{t(n)}g_{n}^{-1}h_{n}^{-1})\]
for $s=(s_{1}, \ldots, s_{n})\in G$. Define a $G\times G$-equivariant map
\[\Psi \colon (G, \pi, \pi_{g})\rightarrow (G, \pi, \pi_{hg})\]
by $\Psi(s)=sh^{-1}$ for $s\in G$. Since $(G, \pi, \tau)$ is also a ME coupling of $G$ and $\Gamma$, we see that $\Phi$ is $G\times \Gamma$-equivariant by Lemma \ref{lem-gen-ME eq-ext}. It is easy to see that $\Psi \circ \Phi(e)=e$ and thus, $\Psi \circ \Phi ={\rm id}$. Therefore, $\tau$ is the restriction of $\pi_{hg}$.
\end{proof}


\section{Lattice embeddings of the mapping class group}

In this final section, we give another application of Corollaries \ref{cor-mcg-ME eq} and \ref{cor-prod-me-eq}, following \cite{furman3}. We describe all lattice embeddings of a finite direct product of mapping class groups into a locally compact second countable (lcsc) group. We fix notations as follows: let $n$ be a positive integer and $M_{i}$ be a surface with  $\kappa(M_{i})>0$ for $i\in \{ 1, \ldots, n\}$. Put $G_{0}=\Gamma(M_{1})^{\diamond}\times \cdots \times \Gamma(M_{n})^{\diamond}$ and $G={\rm Aut}(C(M_{1}))\times \cdots \times {\rm Aut}(C(M_{n}))$. Let $\pi \colon G_{0}\rightarrow G$ be the natural homomorphism.

\begin{thm}\label{thm-lat-emb}
Let $\Gamma$ be a finite index subgroup of $G_{0}$. Suppose that we have a lattice embedding $\sigma \colon \Gamma \rightarrow H$ into a lcsc group $H$. Then there exist the following two maps:
\begin{enumerate}
\item[(i)] an almost $\Gamma \times \Gamma$-equivariant Borel map $\Phi \colon (H, \sigma, \sigma)\rightarrow (G, \pi, \pi)$, which satisfies $\Phi(h_{1}h_{2})=\Phi(h_{1})\Phi(h_{2})$ for a.e.\ $(h_{1}, h_{2})\in H\times H$;
\item[(ii)] a continuous homomorphism $\Phi_{0}\colon H\rightarrow G$ such that $\Phi_{0}(h)=\Phi(h)$ for a.e.\ $h\in H$ and $\Phi_{0}(\sigma(\gamma))=\pi(\gamma)$ for any $\gamma \in \Gamma$. Moreover, $\ker(\Phi_{0})$ is compact.
\end{enumerate}
\end{thm}

\begin{proof}
First, we show the assertion (i). To prove this, we may assume that $M_{i}\neq M_{1, 2}, M_{2, 0}$ for all $i$ by using Lemma \ref{lem-gen-ME eq-ext}. We identify $\Gamma(M_{i})^{\diamond}$ and ${\rm Aut}(C(M_{i}))$ via the natural isomorphism. Applying Corollary \ref{cor-prod-me-eq} to the self ME coupling $H$ of $\Gamma$ with the Haar measure, we can find a bijection $t$ on the set $\{ 1, \ldots, n\}$, an isotopy class $g_{i}$ of a diffeomorphism $M_{t(i)}\rightarrow M_{i}$ and an almost $\Gamma \times \Gamma$-equivariant Borel map 
\[\Phi'\colon (H, \sigma, \sigma)\rightarrow (G, \pi, \pi_{g}).\]
Define a Borel map $F\colon H\times H\rightarrow G$ by
\[F(h_{1}, h_{2})=\Phi'(h_{1}^{-1})^{-1}\Phi'(h_{1}^{-1}h_{2})\Phi'(h_{2})^{-1}\]
for $h_{1}, h_{2}\in H$. Then for any $\gamma \in \Gamma$ and a.e.\ $(h_{1}, h_{2})\in H\times H$, we have
\[F(h_{1}\sigma(\gamma), h_{2})=F(h_{1}, h_{2})=F(h_{1}, h_{2}\sigma(\gamma)^{-1}),\]
\[F(\sigma(\gamma)h_{1}, \sigma(\gamma)h_{2})=\pi_{g}(\gamma)F(h_{1}, h_{2})\pi(\gamma)^{-1}.\]
Thus, $F$ induces a Borel map $f$ from $X=(H/\sigma(\Gamma))\times (H/\sigma(\Gamma))$ to $G$ such that
\[f(\gamma x)=\pi_{g}(\gamma)f(x)\pi(\gamma)^{-1}\]
for any $\gamma \in \Gamma$ and a.e.\ $x\in X$, where the $\Gamma$-action on $X$ is induced from the diagonal one on $H\times H$ through $\sigma$. By projecting the finite $\Gamma$-invariant measure on $X$ to $G$ through $f$, we obtain a finite measure $\mu$ on $G$ invariant under the action of the diagonal subgroup of $\Gamma \times \Gamma$ on $(G, \pi_{g}, \pi)$. It follows that $t={\rm id}$ and $g_{i}\in \Gamma(M_{i})^{\diamond}$. Put $g=(g_{1}, \ldots, g_{n})\in G$. The support of $\mu$ is $\{ g\}$ and $F(h_{1}, h_{2})=g$, that is, $\Phi'(h_{1}^{-1}h_{2})=\Phi'(h_{1}^{-1})g\Phi'(h_{2})$ for a.e.\ $(h_{1}, h_{2})\in H\times H$. 

Define a $G\times G$-equivariant map $\Phi''\colon (G, \pi, \pi_{g})\rightarrow (G, \pi, \pi)$ by $\Phi''(\gamma)=\gamma g$ for $\gamma \in G$. Then the composition $\Phi =\Phi''\circ \Phi'$ is a $\Gamma \times \Gamma$-equivariant Borel map from $(H, \sigma, \sigma)$ to $(G, \pi, \pi)$ and it satisfies $\Phi(h_{1}^{-1}h_{2})=\Phi(h_{1}^{-1})\Phi(h_{2})$ for a.e.\ $(h_{1}, h_{2})\in H\times H$.

Next, we show the assertion (ii). It follows from \cite[Theorems B.2, B.3]{zim2} that there exists a continuous homomorphism $\Phi_{0}\colon H\rightarrow G$ such that $\Phi_{0}(h)=\Phi(h)$ for a.e.\ $h\in H$. For any $\gamma \in \Gamma$ and a.e.\ $h\in H$, we have
\[\pi(\gamma)\Phi(h)=\Phi(\sigma(\gamma)h)=\Phi_{0}(\sigma(\gamma)h)=\Phi_{0}(\sigma(\gamma))\Phi_{0}(h)=\Phi_{0}(\sigma(\gamma))\Phi(h),\]
which implies $\pi(\gamma)=\Phi_{0}(\sigma(\gamma))$ for any $\gamma \in \Gamma$. 

Since $\ker(\Phi_{0})$ is essentially equal to $\Phi^{-1}(e)$, which has finite measure, we see that $\ker(\Phi_{0})$ is compact.
\end{proof}

\noindent {\it Proof of Theorem \ref{thm-int-lattice-emb}.}
It follows from Theorem \ref{thm-lat-emb} that there exists a continuous homomorphism $\Phi_{0}\colon H\rightarrow G$ such that $K=\ker(\Phi_{0})$ is compact. Let $H_{0}=\Phi_{0}^{-1}(\pi(\Gamma))$. These groups satisfy the conditions in the theorem. \hfill $\Box$


\end{document}